\documentclass[10pt,journal]{IEEEtran}
\usepackage{graphicx}          
\usepackage{amsmath}           
\usepackage{amssymb}
\usepackage{algorithm}
\usepackage{algpseudocode}
\usepackage{hyperref}
\usepackage{color}
\usepackage{cite}
\usepackage{epstopdf}
\usepackage[latin1]{inputenc}  

\usepackage{amsthm}

\theoremstyle{definition}

\theoremstyle{remark}
\theoremstyle{plain}

\theoremstyle{remark} \newtheorem{note}{Remark}

\DeclareMathOperator{\E}{E}
\DeclareMathOperator*{\argmin}{arg\,min}

\newcommand{\mbf}[1]{\mathbf{#1}}
\newcommand{\mbs}[1]{\boldsymbol{#1}}

\begin{document}

\title{Online Hyperparameter-Free Sparse \\ Estimation Method}
\author{Dave Zachariah and Petre Stoica\thanks{This work has been partly supported by the Swedish Research Council under contract 621-2014-5874. The authors are with the Department of Information Technology, Uppsala University, Sweden. E-mail: dave.zachariah@it.uu.se and ps@it.uu.se}}

\maketitle

\begin{abstract}
In this paper we derive an online estimator for sparse parameter vectors which, unlike the LASSO approach, does not require the tuning of any hyperparameters. The algorithm is based on a covariance matching approach and is equivalent to a weighted version of the square-root LASSO. The computational complexity of the estimator is of the same order as that of the online versions of regularized least-squares (RLS) and LASSO. We provide a numerical comparison with feasible and infeasible implementations of the LASSO and RLS to illustrate the advantage of the proposed online hyperparameter-free estimator.
\end{abstract}

\section{Introduction}

Estimating a high-dimensional sparse vector of parameters with a few dominant or nonzero
elements has become an important topic in statistics and signal
processing. Applications of sparse estimation include 
spectral analysis \cite{BourguignonEtAl2007_sparsity,KleinEtAl2008_sparsity,StoicaEtAl2011_spicespectral,Stoica&Babu2012_spice}, array
processing \cite{Gorodnitsky&Rao1997_sparse,MalioutovEtAl2005_sparse,StoicaEtAl2011_spicearray},
biomedical analysis \cite{Tibshirani1996_lasso,WuEtAl2009_genome,LuEtAl2011_lasso}, magnetic resonance imaging \cite{LustigEtAl2008_compressedmri,DonevaEtAl2010_compressed_mri},
system identification
\cite{ChenEtAl2009_sparse,KalouptsidisEtAl2011_adaptive,KopsinisEtAl2011_online,Glentis2014_adaptiveslim,ThemelisEtAl2014_variationalbayessparse} and
synthetic aperture radar imaging \cite{Zhu&Bamler2010_tomographic,
  CetinEtAl2014_sparsitysar}.

Many sparse estimation approaches can be
implemented using various computational methods and it is relevant to
formulate estimators that scale well with the size of the
data. Furthermore, in several applications data is obtained as a
stream of measurements and it is desirable to process them
accordingly. Both reasons motivate developing estimation methods that
perform `online' processing, that is, successively refining the estimate of the
sparse parameter vector for each obtained data sample. 
Another common issue with sparse estimation methods is the need for the user to
select or tune critical hyperparameters to strike a balance between
data fidelity and sparsity so as to fit a particular measurement
setup \cite{Malecki&Donoho2010_tunedcompressed,GiraudEtAl_2012_unknownvariance}. This selection is, however, rarely feasible in online
scenarios. Furthermore, when the user has to tune hyperparameters the
outcomes become more arbitrary and the reproducibility of the method
is reduced. Finally, many convex relaxation-based sparse estimation methods are not well
adapted for complex-valued data and parameters and thus they must separate the data into 
real and imaginary parts. This separating approach requires enforcing
pairwise sparsity constraints to avoid performance loss and effectively
doubles the size of computed quantities
\cite{YangEtAl2012_compressedsensingcomplex, MalekiEtAl2013_complexlasso}.

In this paper we develop a sparse
estimation method that addresses the aforementioned
issues. Specifically, 
\begin{itemize}
\item the estimator is implemented online with the
same complexity order as the best existing online methods for sparse estimation. 
\item It automatically adapts to the signal model via
a covariance matching approach and in this way obviates the need
for tuning hyperparameters. 
\item The method can estimate complex-valued parameters as simply as real-valued ones.
\end{itemize}

\emph{Notation:} $\| \cdot \|_1$, $\| \cdot \|_2$ and $\| \cdot \|_F$ denote
the $\ell_1$, $\ell_2$ and Frobenius norms, respectively. Unless otherwise
stated, $\| \cdot \|$ will denote the $\ell_2$-norm and $\| \mbf{x}
\|_{\mbf{W}} = \sqrt{\mbf{x}^* \mbf{W} \mbf{x}}$ where $\mbf{W} \succ
\mbf{0}$ is a positive definite matrix. $[\mbs{\Gamma}]_i$ is
the $i$th column of matrix $\mbs{\Gamma}$ and $\mbs{\Gamma}^\dagger$ is the Moore-Penrose pseudoinverse.

\emph{Abbreviations:} Least squares (LS), regularized least-squares
(\textsc{Rls}), least absolute shrinkage and selector operator (\textsc{Lasso}), sparse
iterative covariance-based estimation (\textsc{Spice}), mean-square error
(MSE), online (\textsc{Ol}).

\section{Background}
 
Consider a sequence of scalar measurements:
\begin{equation}
y_t = \mbf{h}^*_t \mbs{\theta} + w_t \in \mathbb{C}, \quad t = 1, 2,
\dots,
\label{eq:stream}
\end{equation}
where the regressor vector $\mbf{h}_t \in\mathbb{C}^p$ is given, the
unknown sparse parameter vector is $\mbs{\theta} \in
\mathbb{C}^p$ and $w_t$ is zero-mean noise with variance
$\sigma^2$. For the sake of generality we consider complex-valued
variables; any differences that occur in the real-valued case will be addressed below.

Suppose we have obtained $n$ measurements. Then in vector form  we can write
\begin{equation}
\mbf{y}_n = \mbf{H}_n \mbs{\theta} + \mbf{w}_n \in \mathbb{C}^{n},
\label{eq:linearmodel}
\end{equation}
where
\begin{equation*}
\begin{split}
\mbf{H}_n &= \begin{bmatrix} \mbf{h}^*_1 \\ \mbf{h}^*_2 \\ \vdots \\
  \mbf{h}^*_n \end{bmatrix} = \begin{bmatrix} \mbf{c}^{(n)}_1 & \mbf{c}^{(n)}_2 & \cdots &
  \mbf{c}^{(n)}_p \end{bmatrix} \in \mathbb{C}^{n \times p}.
\end{split}
\end{equation*}
To avoid notational clutter we will omit the superindex for the
columns and simply write $\mbf{c}_i$. 
In the following sub-sections we review a few estimators
$\hat{\mbs{\theta}}_n$ of $\mbs{\theta}$ in \eqref{eq:linearmodel}, based on regularizations of the least-squares approach, and
their online formulations that compute $\hat{\mbs{\theta}}_n$ from
$\hat{\mbs{\theta}}_{n-1}$ thus eliminating the need for re-calculating the
estimate from scratch.

\subsection{LS and RLS}

The LS approach  is based on solving the
quadratic problem
\cite{Bjorck1996_numericalls,Soderstrom&Stoica1988_system,KailathEtAl2000_linear}
\begin{equation}
\argmin_{\mbs{\theta} \in \mathbb{C}^p} \; \| \mbf{y}_n - \mbf{H}_n
\mbs{\theta} \|^2_2,
\label{eq:LS}
\end{equation} 
which has the following minimum $\ell_2$-norm solution $\hat{\mbs{\theta}}_n =
\mbf{H}^{\dagger}_n \mbf{y}_n$. If $\mbf{H}_n$ has full column-rank, then
the estimator admits a simple closed-form solution that can be computed by
recursive updates for $n=1, 2, \dots$. This obviates the
need for choosing an initial estimate $\hat{\mbs{\theta}}_0$ or any
hyperparameter (see, e.g. \cite{Stoica&Ahgren2002_exactrls}).

It is more common to consider a regularized LS problem
\begin{equation}
\argmin_{\mbs{\theta} \in \mathbb{C}^p} \; \| \mbf{y}_n - \mbf{H}_n
\mbs{\theta} \|^2_2 + \| \mbs{\theta} \|^2_{\mbs{\Lambda}},
\label{eq:RegLS}
\end{equation} 
with an initial estimate $\hat{\mbs{\theta}}_0 = \mbf{0}$ that
is well-motivated for sparse parameter vectors and with
$\mbs{\Lambda} = \lambda \mbf{I}_p \succeq \mbf{0}$, where $\lambda$ is a hyperparameter chosen
by the user to bias the estimate towards $\mbf{0}$ with the aim of reducing its variance. This estimator
admits an online form $\hat{\mbs{\theta}}_{n} = \hat{\mbs{\theta}}_{n-1} + \mbf{K}_n ( y_n -
\mbf{h}^*_n \hat{\mbs{\theta}}_{n-1})$, where $\mbf{K}_n$ is a matrix
determined from the regressors and $\mbs{\Lambda}$
\cite{KailathEtAl2000_linear}. The computational complexity of this
$\ell_2$-regularized least-squares algorithm is of the order $\mathcal{O}(p^2)$ per sample.

One approach that takes sparsity into account would be to perform
online \textsc{Rls} estimation only on the nonzero components of
$\mbs{\theta}$, if these were known. In \cite{DumitrescuEtAl2012_greedyrls} the
components are successively detected in a greedy manner using
information theoretic criteria at each sample $y_n$. Since the detection
process is subject to errors, the resulting online sparse least-squares
estimate is only an approximation of \eqref{eq:RegLS} applied to the subvector of nonzero coefficients.

\subsection{\textsc{Lasso}}

A substantially different approach than \textsc{Rls} consists of replacing the $\ell_2$-norm regularization
term in \eqref{eq:RegLS} with alternative forms that promote sparsity
directly in the objective itself \cite{Eksioglu&Korhan2011_rlsconvex}. In doing so, sparse solutions can be obtained without the need for concomitantly detecting the nonzero components of $\mbs{\theta}$. This approach to sparse parameter estimation was popularized in
\cite{Tibshirani1996_lasso, ChenEtAl1998_bpdn}. The \textsc{Lasso} estimator solves the following convex problem
\begin{equation}
\argmin_{\mbs{\theta} \in \mathbb{R}^p} \; \| \mbf{y}_n - \mbf{H}_n \mbs{\theta} \|^2_2 + \lambda_n \| \mbs{\theta} \|_1.
\label{eq:lasso}
\end{equation} 
While the solution $\hat{\mbs{\theta}}_n$ does not have a closed-form expression, it can be computed using various numerical methods. Among the
more computationally elegant and scalable methods is the cyclic
minimization strategy of coordinate descent which updates
one element of $\hat{\mbs{\theta}}_n$ at a time in an
iterative manner,
cf. \cite{Fu1998_penalized,FriedmanEtAl2007_pathwise} and references therein.

One way of formulating an online solution is to interpret
\eqref{eq:lasso} as a penalized maximum likelihood estimator, assuming 
Gaussian noise in \eqref{eq:stream}. Then it is possible to formulate an iterative expectation
maximization algorithm with recursively updated quantities using
auxiliarly variables \cite{BabadiEtAl2010_sparls}. The drawback,
however, is that an additional hyperparameter, besides $\lambda_n$ in \eqref{eq:lasso}, needs to be tuned. Another way of dynamically updating the estimate $\hat{\mbs{\theta}}_{n-1}$
from $\hat{\mbs{\theta}}_n$ is the method of homotopy
\cite{Garrigues&Ghaoi2009_homotopylasso,Asif&Romberg2010_homotopylasso},
whereby the cost function in \eqref{eq:lasso} with a fixed $\lambda_n
\equiv \lambda$ is modified into $\| \mbf{y}_{n-1} - \mbf{H}_{n-1}
\mbs{\theta} \|^2_2 + \epsilon | y_n - \mbf{h}^*_n \mbs{\theta}|^2 +
\lambda \| \mbs{\theta} \|_1$. As the scalar parameter $\epsilon \in
[0,1]$ is varied from 0 to 1, the transition from $\hat{\mbs{\theta}}_{n-1}$ to $\hat{\mbs{\theta}}_n$ can be
computed more efficiently than recalculating $\hat{\mbs{\theta}}_n$
from scratch thereby enabling an online formulation.

For the real-valued case, an elegant online formulation is found in
\cite{AngelosanteEtAl2010_onlinelasso}, which is based on the cyclic
minimization strategy mentioned above. The cost function in
\eqref{eq:lasso} can be written equivalently as $\mbs{\theta}^\top
\mbs{\Gamma}^n \mbs{\theta} - 2 \mbs{\theta}^\top \mbs{\rho}^n +
\lambda_n \| \mbs{\theta} \|_1$, ignoring any constant, where
$\mbs{\Gamma}^n = \mbf{H}^*_n \mbf{H}_n$ and $\mbs{\rho}^n =
\mbf{H}^*_n \mbf{y}_n$ can be computed recursively. Then, starting
from an initial estimate $\hat{\mbs{\theta}}_0$, the elements of
$\hat{\mbs{\theta}}_n$ are updated for each sample by solving
\begin{equation*}
\hat{\theta}_i = \argmin_{\theta_i} \; \Gamma^n_{ii} \theta^2_i - 2\tilde{\rho}^n_i \theta_i + \lambda_n |\theta_i|
\end{equation*}
in closed form for $i=1,\dots,p$, where $\tilde{\rho}^n_i = \rho^n_i -
\sum_{j\neq i} \Gamma^n_{ij} \check{\theta}_j$ and $\check{\theta}_j$
denotes the current estimate. The complexity of the full online cyclic minimization \textsc{Lasso} is $\mathcal{O}(p^2)$ per sample.

Under certain conditions on the regressors, sparsity of $\mbs{\theta}$,
and noise, it is possible to prove that the \textsc{Lasso} estimator
possesses `oracle' properties. That is, asymptotically it can identify
the support set of $\mbs{\theta}$ and perform as well as \textsc{Rls}
applied to the nonzero coefficients of the parameter vector, cf. \cite{Fuchs2005_recovery,DonohoEtAl2006_stable,Tropp2006_justrelax,AngelosanteEtAl2010_onlinelasso}. This,
however, requires selecting the hyperparameter $\lambda_n$ based on the knowledge of the noise variance $\sigma^2$ which is rarely feasible in practical (online) scenarios.

\subsection{Square-root \textsc{Lasso}}

To circumvent the need to know $\sigma$ in the \textsc{Lasso}, a subtle modification of \eqref{eq:lasso} was proposed in \cite{BelloniEtAl2011_squarerootlasso},
\begin{equation}
\argmin_{\mbs{\theta} \in \mathbb{R}^p} \; \| \mbf{y}_n - \mbf{H}_n
\mbs{\theta} \|_2 + \lambda_n \| \mbs{\theta} \|_1,
\label{eq:sqrtlasso}
\end{equation} 
where the first term, containing the residuals, is the square-root of that in
\eqref{eq:lasso}. As argued in \cite{BelloniEtAl2011_squarerootlasso},
near-oracle performance for both \eqref{eq:lasso} and \eqref{eq:sqrtlasso} can be achieved when $\lambda_n$ is chosen as the smallest
value that dominates the gradient of the first term, when
evaluated at the true $\mbs{\theta}$. At this point, the gradient
captures the estimation errors arising from noise alone. However, by
re-parameterizing \eqref{eq:linearmodel} as $\mbf{y}_n = \mbf{H}_n
\mbs{\theta} + \sigma \mbs{\varepsilon}$, where
$\E[\mbs{\varepsilon}\mbs{\varepsilon}^\top] = \mbf{I}_n$, it is seen
that the gradients of the first terms in \eqref{eq:lasso} and
\eqref{eq:sqrtlasso} differ in one crucial respect; namely the latter
does not depend on $\sigma$ thus rendering the choice of $\lambda_n$ for
\eqref{eq:sqrtlasso} invariant to the noise level. 

Another way to address the dependence on $\sigma$ is to estimate it
\cite{Sun&Zhang2012_scaledsparse}. The square-root \textsc{Lasso} estimator in
\eqref{eq:sqrtlasso} can in fact be interpreted as an
$\ell_1$-penalized joint estimator of $\mbs{\theta}$ and $\sigma$
used in robust regression. Suppose $\gamma(\cdot)$ is a convex loss function
of the normalized residuals $(y_t -
\mbf{h}^\top_t\mbs{\theta}) /\sigma$. Then the concomitant M-estimator of location and
scale, $\mbs{\theta}$ and $\sigma$, is given by 
\cite[ch.7]{Huber2011_robust}
 \begin{equation}
 \argmin_{\mbs{\theta}, \: \sigma } \; \frac{1}{n} \sum^n_{t=1} \left[
  \gamma\left( \frac{y_t - \mbf{h}^\top_t\mbs{\theta}}{\sigma} \right) + 
   a \right] \sigma ,
 \label{eq:sqrtlasso_huber}
 \end{equation}
where $a > 0$ is a user-defined parameter. In robust regression, various loss
functions are considered to mitigate noise outliers. For a squared-error loss $\gamma(x) =
x^2$, we obtain the minimizer $\hat{\sigma} = \|
\mbf{y}_n - \mbf{H}_n \mbs{\theta} \|_2 / \sqrt{na}$ in closed form. Penalizing the
M-estimator in \eqref{eq:sqrtlasso_huber} by $ \lambda_n \|
\mbs{\theta} \|_1$ and concentrating out the minimizing $\sigma$
with $a=n/4$ yields \eqref{eq:sqrtlasso}.

While an efficient choice of $\lambda_n$ in \eqref{eq:sqrtlasso} is
independent of $\sigma$, the user input is still required; furthermore, the choices of $\lambda_n$ in
\cite{BelloniEtAl2011_squarerootlasso} are predicated on the assumption that each column
of $\mbf{H}_n$ has unit norm. Such a rescaling of the regressors may not be practical in
an online scenario. Note that a cyclic minimization algorithm for the
convex square-root \textsc{Lasso} has been presented in the
supplementary material of \cite{BelloniEtAl2011_squarerootlasso} (albeit only for the real-valued case and without any derivation) but an online implementation has not yet been formulated.

\subsection{\textsc{Spice} as weighted square-root \textsc{Lasso}}

Let us now consider the estimation problem from a statistical point of
view. Suppose $\mbs{\theta}$ is a zero-mean random variable with
covariance matrix $\mbf{P} \succ \mbf{0}$. Then the linear estimator that minimizes the mean square error $\E_{y,\theta}[\|\mbs{\theta} - \hat{\mbs{\theta}}_n \|^2_2]$ is obtained by solving
\begin{equation}
\begin{split}
\argmin_{\mbs{\theta}} \; \frac{1}{\sigma^2} \| \mbf{y}_n -
\mbf{H}_n \mbs{\theta} \|^2_2 + \| \mbs{\theta} \|^2_{\mbf{P}^{-1}},
\end{split}
\label{eq:lmmse_problem}
\end{equation}
and can be written in closed form as \cite{vanTrees2013_detection,Soderstrom&Stoica1988_system,KailathEtAl2000_linear}
\begin{equation}
\begin{split}
\hat{\mbs{\theta}}_n &= \mbf{P}\mbf{H}^*_n( \mbf{H}_n\mbf{P}\mbf{H}^*_n + \sigma^2 \mbf{I}_n )^{-1} \mbf{y}_n \\
&= ( \mbf{H}^*_n \mbf{H}_n + \sigma^2 \mbf{P}^{-1} )^{-1} \mbf{H}^*_n \mbf{y}_n .
\end{split}
\label{eq:lmmse}
\end{equation}
In the problem under consideration, however, neither $\mbf{P}$ nor
$\sigma^2$ is known. By treating them as unknown parameters, they can
be estimated by a covariance-matching approach (e.g.,
\cite{Anderson1989_linear,OtterstenEtAl1998_covariance}) and then used in
\eqref{eq:lmmse}. 

For reasons of parsimony and tractability we do not model any correlations
between the elements of $\mbs{\theta}$ and hence $\mbf{P}$ is a $p
\times p$
diagonal matrix. Now consider the covariance matrix of the data $\mbf{R}_n = \E[\mbf{y}_n \mbf{y}^*_n] =  \mbf{H}_n \mbf{P} \mbf{H}^*_n +
\sigma^2 \mbf{I}_n$, which is a function of $\mbf{P}$ and
$\sigma^2$. We choose these $p+1$ nonnegative parameters to match the
covariance of the observed data, by minimizing the criterion $$\|
\mbf{R}^{-1/2}_n(\mbf{y}_n\mbf{y}^*_n - \mbf{R}_n)  \|^2_F,$$ with
respect to $\mbf{P}$ and $\sigma^2$. This criterion is
the basis of the sparse iterative covariance-based estimation (\textsc{Spice}) framework.

Using this covariance-matching approach is equivalent to solving for the parameters
jointly in the following augmented problem
\begin{equation}
\begin{split}
\argmin_{\mbs{\theta}, \: \mbf{P}, \: \sigma^2} \; &\frac{1}{\sigma^2} \| \mbf{y}_n -
\mbf{H}_n \mbs{\theta}\|^2_2 + \| \mbs{\theta} \|^2_{\mbf{P}^{-1}} \\
&\quad + \text{tr}\left\{ \mbf{H}_n
  \mbf{P} \mbf{H}^*_n  + \sigma^2 \mbf{I}_n\right\},
\end{split}
\label{eq:lmmseaugprob}
\end{equation}
which is similar in form to \eqref{eq:lmmse_problem} but contains the additional
term $\text{tr} \{ \mbf{R}_n  \} = \text{tr}\left\{ \mbf{H}_n
  \mbf{P} \mbf{H}^*_n  + \sigma^2 \mbf{I}_n\right\}$. (See Appendix~A for a proof of this equivalence.) Furthermore, following
\cite{RojasEtAl2013_spicenote,Babu&Stoica2014_connection,StoicaEtAl2014_weightedspice}
it can be shown that solving for $\mbf{P}$ and $\sigma^2$, and concentrating them out from \eqref{eq:lmmseaugprob}, results in
\begin{equation}
\argmin_{\mbs{\theta} \in \mathbb{C}^p} \; \| \mbf{y}_n - \mbf{H}_n
\mbs{\theta} \|_2 +  \| \mbf{D}_n \mbs{\theta} \|_1,
\label{eq:weightedsqrtLASSO}
\end{equation} 
where
\begin{equation*}
\mbf{D}_n = 
\text{diag}\left( \sqrt{\frac{\| \mbf{c}_1 \|^2_2}{n}}, \dots, \sqrt{\frac{\| \mbf{c}_p \|^2_2}{n}} \right).
\end{equation*}
Eq. \eqref{eq:weightedsqrtLASSO} can be interpreted as a weighted, hyperparameter-free
square-root \textsc{Lasso}. As is the case with the square-root \textsc{Lasso}, online formulations of \eqref{eq:weightedsqrtLASSO} have not appeared in the literature.

\subsection{Problem formulation}

We have reviewed several approaches
to sparse parameter estimation as well as some of their
interconnections and limitations. Note that all of the estimators
considered above involve convex minimization problems. The
$\ell_2$ and $\ell_1$-penalized forms of \eqref{eq:LS} in
\eqref{eq:RegLS} and \eqref{eq:lasso} have concise online formulations
but require the careful selection of hyperparameters. Furthermore, an
efficient choice depends on the unknown noise power $\sigma^2$.  The hyperparameters choice is rendered invariant to $\sigma^2$ by the change in
\eqref{eq:sqrtlasso}. Moreover, this selection
is entirely avoided in \eqref{eq:weightedsqrtLASSO} using the \textsc{Spice} approach.

The goal of the remainder of the
paper is to formulate an online \textsc{Spice} estimator for the sparse vector
$\mbs{\theta}$ (see \eqref{eq:weightedsqrtLASSO}) given data $\{ y_t, \mbf{h}_t \}^{n}_{t=1}$. This estimator, denoted `\textsc{Ol-Spice}', obviates the need for user-defined hyperparameters, treats the
complex-valued case as simply as the real-valued one, and is of the same complexity order as the online solutions of \eqref{eq:RegLS} and \eqref{eq:lasso}. In the numerical example section we provide results comparing the aforementioned online estimators, viz. \textsc{Ol-Rls}, \textsc{Ol-Lasso} and \textsc{Ol-Spice}.

\section{Online \textsc{Spice}}
\label{sec:onlinespice}

First we
formulate a low-complexity cyclic minimization algorithm for the cost
function in \eqref{eq:weightedsqrtLASSO}. Then, using this result we
derive an online estimator which sequentially processes a stream of data with complexity $\mathcal{O}(p^2)$ per sample.

\subsection{Cyclic minimization}

Let the cost function in \eqref{eq:weightedsqrtLASSO} be denoted as $J(\mbs{\theta}) =  \| \mbf{y}_n - \mbf{H}_n \mbs{\theta} \|_2 +
 \| \mbf{D}_n \mbs{\theta} \|_1$. We begin by minimizing
$J(\mbs{\theta})$ with respect to one component $\theta_i$ at a time. Let
$\tilde{\mbf{y}}_i \triangleq \mbf{y} - \sum_{k \neq i} \mbf{c}_k
\theta_k $ (omitting the index $n$ to lighten the notation); then the cost function can be re-written as
\begin{equation}
J(\theta_i) = (\| \tilde{\mbf{y}}_i- \mbf{c}_i \theta_i \|^2_2)^{1/2} + d_{ii}
|\theta_i| + K,
\label{eq:globalcost_alt}
\end{equation}
where $d_{ii}$ is the $i$th diagonal element of $\mbf{D}_n$ and
$K=\sum_{k \neq i} d_{kk}|\theta_k|$ is a constant. To tackle this
scalar minimization problem we reparameterize the $i$th variable in
polar form $\theta_i = r_i e^{j \varphi_i}$ where $r_i \geq 0$ and $\varphi_i \in [-\pi, \pi )$
(or $\varphi_i \in \{0, \pi \}$ when $\theta_i$ is real-valued). This
enables the following reformulation of the quadratic term in \eqref{eq:globalcost_alt}:
\begin{equation}
\begin{split}
\| \tilde{\mbf{y}}_i- \mbf{c}_i \theta_i \|^2_2 &= \| \tilde{\mbf{y}}_i - \mbf{c}_i r_i e^{j \varphi_i} \|^2_2 \\
&= \| \tilde{\mbf{y}}_i \|^2_2 + \| \mbf{c}_i r_i   e^{j \varphi_i} \|^2_2 - 2
\text{Re}\{ r_i \mbf{c}^*_i \tilde{\mbf{y}}_i e^{-j\varphi_i} \}  \\
&= \| \tilde{\mbf{y}}_i \|^2_2 + \| \mbf{c}_i \|^2_2 r^2_i  \\
&\quad - 2 r_i |\mbf{c}^*_i \tilde{\mbf{y}}_i|
\cos( \text{arg}(\mbf{c}^*_i \tilde{\mbf{y}}_i) - \varphi_i ) .
\end{split}
\label{eq:quadraticterm}
\end{equation}
Inserting \eqref{eq:quadraticterm} into $J (\theta_i)$ and noting that
$|\theta_i| = |r_i e^{j \varphi_i}| = r_i$, we obtain the following criterion as a function of
$r_i$ and $\varphi_i$,
\begin{equation}
\begin{split}
& J(r_i, \varphi_i) \\
&= \left(\| \tilde{\mbf{y}}_i \|^2_2 + \| \mbf{c}_i \|^2_2 r^2_i  - 2 r_i |\mbf{c}^*_i \tilde{\mbf{y}}_i|
\cos( \text{arg}(\mbf{c}^*_i \tilde{\mbf{y}}_i) - \varphi_i )
\right)^{1/2} \\
&\quad + d_{ii} r_i .
\label{eq:J_r_phi}
\end{split}
\end{equation}
The minimizing $\varphi_i$ is simply
\begin{equation}
\hat{\varphi}_i =  \text{arg}(\mbf{c}^*_i \tilde{\mbf{y}}_i) ,
\label{eq:phi_hat}
\end{equation}
whether the data is complex or real-valued.

Next, let
\begin{equation}
\begin{split}
\alpha_i &\triangleq \| \tilde{\mbf{y}}_i \|^2\\
\beta_i &\triangleq \| \mbf{c}_i \|^2\\
\gamma_i &\triangleq  |\mbf{c}^*_i \tilde{\mbf{y}}_i|,
\end{split}
\label{eq:updatevariables}
\end{equation}
so that we can write \eqref{eq:J_r_phi} as 
\begin{equation}
\begin{split}
J(r_i, \hat{\varphi}_i) &= \left( \alpha_i + \beta_i r^2_i  - 2 \gamma_i r_i \right)^{1/2} + d_{ii}
r_i .
\end{split}
\label{eq:J_r_phi_alt}
\end{equation}
Note that by the Cauchy-Schwarz inequality 
\begin{equation}
\alpha_i \beta_i - \gamma^2_i \geq 0.
\label{eq:CauchySchwarz}
\end{equation}
Equality in \eqref{eq:CauchySchwarz} occurs only when
$\tilde{\mbf{y}}_i$ is colinear with $\mbf{c}_i$. Inserting
\eqref{eq:linearmodel} into $\tilde{\mbf{y}}_i$ one obtains $\tilde{\mbf{y}}_i =
 \sum_{k\neq i} \mbf{c}_{k} \tilde{\theta}_k + \mbf{c}_i \theta_i +
 \mbf{w} $, where $\tilde{\theta}_k$ denote estimation errors when
 holding the remaining coefficients constant. Due to the random noise
 $\mbf{w}$, $\tilde{\mbf{y}}_i$ and  $\mbf{c}_i$ will not be colinear, making the inequality \eqref{eq:CauchySchwarz} strict, with probability 1.

We now show that \eqref{eq:J_r_phi_alt} is convex and derive the minimizing $r\geq 0$ of
this function (dropping the index
$i$, in what follows, for notational simplicity). The first-order derivative is
\begin{equation}
\begin{split}
\frac{dJ}{dr} &= \frac{\beta r - \gamma}{\left(  \beta
    r^2  - 2 \gamma r + \alpha \right)^{1/2}} + d,
\end{split}
\label{eq:firstderivative}
\end{equation}
where the quadratic expression in the denominator can be factored as
\begin{equation}
\begin{split}
\beta r^2  - 2 \gamma r + \alpha &= \beta \left[ \left(r -
  \frac{\gamma}{\beta} \right)^2 + \left(\frac{\alpha}{\beta} -
  \frac{\gamma^2}{\beta^2}\right) \right].
\end{split}
\label{eq:quadraticdenom}
\end{equation}
Given the strict inequality in \eqref{eq:CauchySchwarz} it follows that the right-hand side of
\eqref{eq:quadraticdenom}, and therefore the denominator of  \eqref{eq:firstderivative}, is strictly
positive. 
The second-order derivative can be expressed as
\begin{equation*}
\begin{split}
\frac{d^2J}{dr^2} &= \frac{\beta}{\left(  \beta r^2  - 2 \gamma r + \alpha
  \right)^{1/2}} - \frac{(\beta r - \gamma)^2}{\left(  \beta
    r^2  - 2 \gamma r + \alpha \right)^{3/2}} \\
&= \frac{1}{\left(  \beta
    r^2  - 2 \gamma r + \alpha \right)^{3/2}} \left( \beta \left(  \beta
    r^2  - 2 \gamma r + \alpha \right) - (\beta r -\gamma)^2 \right ) \\
&= \frac{1}{\left(  \beta   r^2  - 2 \gamma r + \alpha \right)^{3/2}}
\left( \alpha \beta - \gamma^2 \right).
\\
\end{split}
\end{equation*}
Note that the above equation is positive, in view of \eqref{eq:CauchySchwarz}. Thus the function
\eqref{eq:J_r_phi_alt} is convex and the minimizer $r > 0$ is given by
its stationary point; or else $r=0$. 

The stationary point, for which $dJ/dr =
0$, can be found by solving (see \eqref{eq:firstderivative}):
\begin{equation*}
\begin{split}
(\beta r - \gamma) &= -d (\beta r^2 - 2\gamma r + \alpha)^{1/2} ,
\end{split}
\end{equation*}
which leads to the condition $0 \leq r \leq \frac{\gamma}{\beta}$ given that both
sides must be negative. By squaring both sides of the above expression we can write 
\begin{equation*}
\beta^2 \left( r - \frac{\gamma}{\beta} \right)^2 = d^2 \beta
\left[ \left(r - \frac{\gamma}{\beta} \right)^2  +
  \left(\frac{\alpha}{\beta} - \frac{\gamma^2}{\beta^2} \right) \right],
\end{equation*}
using \eqref{eq:quadraticdenom},
or
\begin{equation*}
\left( r - \frac{\gamma}{\beta} \right)^2 = \frac{d^2}{\beta - d^2}
\left( \frac{\alpha \beta - \gamma^2}{\beta^2} \right),
\end{equation*}
where $d^2 / (\beta - d^2) = (\|\mbf{c}\|^2/n)/(\|\mbf{c} \|^2 -
\|\mbf{c}\|^2/n) = 1/(n-1)$. Therefore for $n>1$, we can write the solution more compactly as (reinstating the dependence on $i$):
\begin{equation}
\begin{split}
\hat{r}_i 
&= \frac{\gamma_i}{\beta_i} - \frac{1}{\beta_i} \left( \frac{ \alpha_i \beta_i - \gamma^2_i}{n-1} \right)^{1/2},
\end{split} 
\label{eq:r_hat}
\end{equation}
given the fact that $0 \leq r_i \leq \frac{\gamma_i}{\beta_i}$.

Finally, we summarize the element-wise minimizer of \eqref{eq:globalcost_alt} as
\begin{equation}
\hat{\theta}_i = 
\begin{cases}
\hat{r}_i e^{j \hat{\varphi}_i}, & \text{if }\sqrt{n-1} \gamma_i >
\sqrt{ \alpha_i \beta_i - \gamma^2_i }\\
0 , & \text{else}
\end{cases}
\label{eq:theta_hat_i}
\end{equation}
using \eqref{eq:r_hat} and \eqref{eq:phi_hat}. Updating each element $\hat{\theta}_i$ while holding
the remaining elements constant will monotonically reduce the convex cost function
\eqref{eq:globalcost_alt}. Thus
repeating \eqref{eq:theta_hat_i} for $i=1, \dots, p$ results in a computationally simple cyclic minimizer.

\subsection{Online formulation}

We now derive a method for efficiently computing
\eqref{eq:theta_hat_i}, given the current estimate which we denote by
$\check{\mbs{\theta}}$, at any
$n$, for notational simplicity. At $n=0$,
the estimate is initialized as $\check{\mbs{\theta}} = \mbf{0}$.
We note that the variables in \eqref{eq:phi_hat} and \eqref{eq:updatevariables} depend on quantities whose dimensions grow
with $n$; namely, $\tilde{\mbf{y}}_i$ and $\mbf{c}_i$. By introducing
recursively computed variables we derive an estimate update that keeps the
complexity and memory storage constant at each sample and is of the
same complexity order as online \textsc{Rls} and \textsc{Lasso}.

First, we introduce the auxiliary variable $\mbf{z}_n = \mbf{y}_n - \mbf{H}_n
\check{\mbs{\theta}}$, which will subsequently be eliminated as we
proceed in the derivation. Then we
can write the following identity $\tilde{\mbf{y}}_i = \mbf{z}_n + \mbf{c}_i
\check{\theta}_i$, which enables the variables in \eqref{eq:updatevariables}
to be expressed as
\begin{equation}
\begin{split}
\alpha_i &= \| \tilde{\mbf{y}}_i \|^2 \\
&= \| \mbf{z}_n + \mbf{c}_i \check{\theta}_i \|^2 \\
&= \| \mbf{z}_n \|^2 + \| \mbf{c}_i \|^2 |\check{\theta}_i|^2 + 2 \text{Re}\{
\check{\theta}^*_i \mbf{c}^*_i \mbf{z}_n \}  \\
\beta_i &= \| \mbf{c}_i \|^2 \\
\gamma_i &= |\mbf{c}^*_i \tilde{\mbf{y}}_i| \\
&= |\mbf{c}^*_i (\mbf{z}_n + \mbf{c}_i \check{\theta}_i ) |.
\end{split}
\label{eq:updatevariables_z}
\end{equation}
Next, introduce the auxiliary variables
\begin{equation}
\begin{split}
\eta_n &\triangleq \| \mbf{z}_n \|^2 \\
\mbs{\zeta}_n &\triangleq \mbf{H}^*_n \mbf{z}_n 
\end{split}
\label{eq:auxiliaryvariables}
\end{equation}
and the recursively computed variables
\begin{equation}
\begin{split}
\mbs{\Gamma}^n &\triangleq \mbf{H}^*_n \mbf{H}_n = \mbs{\Gamma}^{n-1} + \mbf{h}_n \mbf{h}^*_n \\
\mbs{\rho}^n &\triangleq \mbf{H}^*_n \mbf{y}_n = \mbs{\rho}^{n-1} + \mbf{h}_n {y}_n \\
\kappa^n &\triangleq \mbf{y}^*_n \mbf{y}_n = \kappa^{n-1} + |y_n|^2,
\end{split}
\label{eq:recursivevariables}
\end{equation}
that are initialized as $\mbf{0}$. 
Then \eqref{eq:updatevariables_z} can be simplified as follows:
\begin{equation}
\begin{split}
\alpha_i &=\eta_n + \Gamma^n_{ii}|\check{\theta}_i|^2 + 2\text{Re}\{
\check{\theta}^*_i \zeta_i \}  \\
\beta_i &= \Gamma^n_{ii}\\
\gamma_i &= |\zeta_i + \Gamma^n_{ii} \check{\theta}_i | ,
\end{split}
\label{eq:updatevariables_compact}
\end{equation}
where $\zeta_i$ denotes the $i$th element of $\mbs{\zeta}_n$.
Similarly, \eqref{eq:phi_hat} can be expressed as
\begin{equation}
\hat{\varphi}_i = \text{arg}(\zeta_i + \Gamma^n_{ii} \check{\theta}_i).
\label{eq:phi_hat_update}
\end{equation}
Therefore the computation of \eqref{eq:theta_hat_i} can be
expressed in terms of  \eqref{eq:auxiliaryvariables}, \eqref{eq:recursivevariables} and the current estimate $\check{\theta}_i$.

Once $\hat{\theta}_i$ has been computed, the current estimate must be updated along with the auxiliary variables to compute the subsequent coefficients of $\hat{\mbs{\theta}}$. The variable
$\mbf{z}_n$ can easily be updated as $\mbf{z}'_n = \mbf{z}_n + \mbf{c}_i
(\check{\theta}_i - \hat{\theta}_i)$, and it follows that the update
of \eqref{eq:auxiliaryvariables} equals
\begin{equation*}
\begin{split}
\eta'_n &= \| \mbf{z}'_n \|^2 \\
&= \eta_n  + {\Gamma}^n_{ii} |\check{\theta}_i - \hat{\theta}_i|^2 + 2
\text{Re}\{  (\check{\theta}_i - \hat{\theta}_i)^* \zeta_i \}\\
\mbs{\zeta}'_n &= \mbf{H}^*_n \mbf{z}'_n \\
&= \mbs{\zeta}_n + [\mbs{\Gamma}^n]_i (\check{\theta}_i - \hat{\theta}_i),
\end{split}
\end{equation*}
which involves a small number of scalar operations and an addition
of two $p \times 1$ vectors. The variable $\mbf{z}_n$ can now be
eliminated, initializing the auxiliary variables
\eqref{eq:auxiliaryvariables} for each sample $n$ as $\eta_n =
\kappa^n + \check{\mbs{\theta}}^* \mbs{\Gamma}^n \check{\mbs{\theta}} - 2 \text{Re}\{
    \check{\mbs{\theta}}^* \mbs{\rho}^n \}$ and $\mbs{\zeta}^n = \mbs{\rho}^n
    -  \mbs{\Gamma}^n \check{\mbs{\theta}}$. We summarize \textsc{Ol-Spice} in Algorithm \ref{alg:onlinesparse}. The algorithm specifies the update
of the estimate for each new sample $y_n$ and is initialized at $n=0$ by $\check{\mbs{\theta}} = \mbf{0}$. The cyclic computation of
$\hat{\theta}_i$, $i = 1, \dots, p$ is terminated after $L \geq 1$ repetitions
per sample, cf. line~14 in Algorithm \ref{alg:onlinesparse}. 

In sum, by introducing the auxiliary variables we can
maintain constant storage and a computational complexity of order
$\mathcal{O}(Lp^2)$ per sample. Since $L\geq 1$ is a constant independent of
$p$, this is the same complexity order as online \textsc{Rls}
and \textsc{Lasso}. As reported below, $L=1$ performs well in practice.
Other update strategies, akin to those
 considered in \cite{AngelosanteEtAl2010_onlinelasso}, can be explored
 in online applications where complexity needs to be further reduced.

\begin{algorithm}
  \caption{: \textbf{Online \textsc{Spice}}} \label{alg:onlinesparse}
\begin{algorithmic}[1]
    \State Input: $y_n$, $\mbf{h}_n$ and $\check{\mbs{\theta}}$
    \State $\mbs{\Gamma} := \mbs{\Gamma} + \mbf{h}_n \mbf{h}^*_n$
    \State $\mbs{\rho} := \mbs{\rho} + \mbf{h}_n y_n$ 
    \State $\kappa := \kappa + |y_n|^2$
    \State $\eta = \kappa + \check{\mbs{\theta}}^* \mbs{\Gamma} \check{\mbs{\theta}} - 2 \text{Re}\{
    \check{\mbs{\theta}}^* \mbs{\rho} \}$   
  \State $\mbs{\zeta} = \mbs{\rho} - \mbs{\Gamma} \check{\mbs{\theta}} $
    \Repeat
        \State $i = 1, \dots, p$
        \State Compute \eqref{eq:updatevariables_compact} and \eqref{eq:phi_hat_update}
        \State Compute $\hat{\theta}_i$ using \eqref{eq:theta_hat_i}
        \State $\eta := \eta  + \Gamma_{ii} |\check{\theta}_i - \hat{\theta}_i|^2 + 2
\text{Re}\{  (\check{\theta}_i - \hat{\theta}_i)^* \zeta_i \}$        
        \State  $\mbs{\zeta}  := \mbs{\zeta} + [\mbs{\Gamma}]_i (\check{\theta}_i - \hat{\theta}_i)$ 
        \State $\check{\theta}_i := \hat{\theta}_i$
    \Until{ number of iterations equals $L$ }
    \State Output: $\hat{\mbs{\theta}}$
\end{algorithmic}
\end{algorithm}

\begin{note}
At any $n$, the output of the algorithm converges to the global minimizer  \eqref{eq:weightedsqrtLASSO} as $L \rightarrow \infty$ which follows from the above analysis of the convex minimization problem. A convergence analysis for finite $L$ and $n \rightarrow
  \infty$ is, however, nontrivial, cf. \cite[ch.~9]{Soderstrom&Stoica1988_system}. 
\end{note}

\begin{note}
The original
\textsc{Spice} batch algorithms
\cite{StoicaEtAl2011_spicespectral,StoicaEtAl2011_spicearray}, with
uniform noise variance, and the above online formulation solve the same
convex problem iteratively. The former uses an initial batch estimate
whereas the latter is initialized by setting $\check{\mbs{\theta}}=\mbf{0}$. A more important difference, however, is that
the former requires repeated inversions of $n \times n$ matrices, each
of which is of complexity $\mathcal{O}(n^3)$, whereas the latter
requires none. This renders batch \textsc{Spice} intractable when $n$ takes on large values (such as $n>1000$ for a regular PC) and
precludes its use in scenarios considered in this work.
\end{note}

\begin{note}
We note that the approach employed to derive
\textsc{Ol-Spice} also enables an alternative formulation of
\textsc{Ol-Lasso} that treats the complex-valued case as simply as the
real-valued one. See Appendix~B for a derivation.
\end{note}

\section{Numerical evaluation}

In this section we compare the derived \textsc{Ol-Spice} with
feasible and infeasible versions of the \textsc{Ol-Rls} and \textsc{Ol-Lasso} \cite{AngelosanteEtAl2010_onlinelasso}.

The infeasible \textsc{Ol-Rls} is implemented by processing only the (unknown)
subset of nonzero coefficients in $\mbs{\theta}$, whereas the feasible
\textsc{Ol-Rls} processes the entire vector, with the regularization parameter set arbitrarily to
$\lambda = 1$. The infeasible \textsc{Ol-Lasso} is implemented by setting
$\lambda_n = \sqrt{2 \sigma^2 n \log p}$, which is proportional to the
(unknown) noise level \cite{AngelosanteEtAl2010_onlinelasso}, whereas
for the feasible \textsc{Ol-Lasso} we set $\lambda_n = \sqrt{n \log p}$.

The performance of the estimators was evaluated using the
normalized mean-square error
\begin{equation}\label{eq:NMSE}
\text{NMSE} \triangleq \frac{\E_{y, \theta}[\| \mbs{\theta} -
  \hat{\mbs{\theta}} \|^2_2]}{\E_{\theta}[\| \mbs{\theta} \|^2_2]}.
\end{equation}
When $\mbs{\theta}$ is an unknown deterministic variable, the expectation with respect to it in \eqref{eq:NMSE} should be be removed. Note that an NMSE value below 0~dB quantifies the error reduction from the initial guess $\hat{\mbs{\theta}} = \mbf{0}$. The NMSE was evaluated using
100~Monte Carlo simulations. We used a PC with Intel i7 3.4~GHz CPU and
16~GB RAM. The algorithms were implemented in
\textsc{Matlab} without any special code optimization or hardware acceleration.

\emph{Remark:} In the interest of reproducible research we have made
the codes for \textsc{Ol-Spice}, as well as for the presented numerical
experiments, available at \textcolor{black}{\texttt{https://www.it.uu.se/katalog/davza513}}.

\subsection{Real-valued example: random regressors}

To illustrate the performance of the estimators we consider a scenario
with the real-valued regressor elements $\mbf{h}_t$ in \eqref{eq:stream} drawn from identical and
independent Gaussian distributions (i.i.d.) with zero mean and unit variance.
The signal to noise ratio is defined as
\begin{equation*}
\text{SNR}  = \frac{\min_{i \in \mathcal{S}} \: \E_{\theta}[|\theta_i|^2]}{\sigma^2},
\end{equation*}
where $\mathcal{S}$ is the support set of $\mbs{\theta}  \in \mathbb{R}^p$ and $p=500$. We use Gaussian noise throughout all experiments.

\begin{figure*}
  \begin{center}
    \includegraphics[width=1.95\columnwidth]{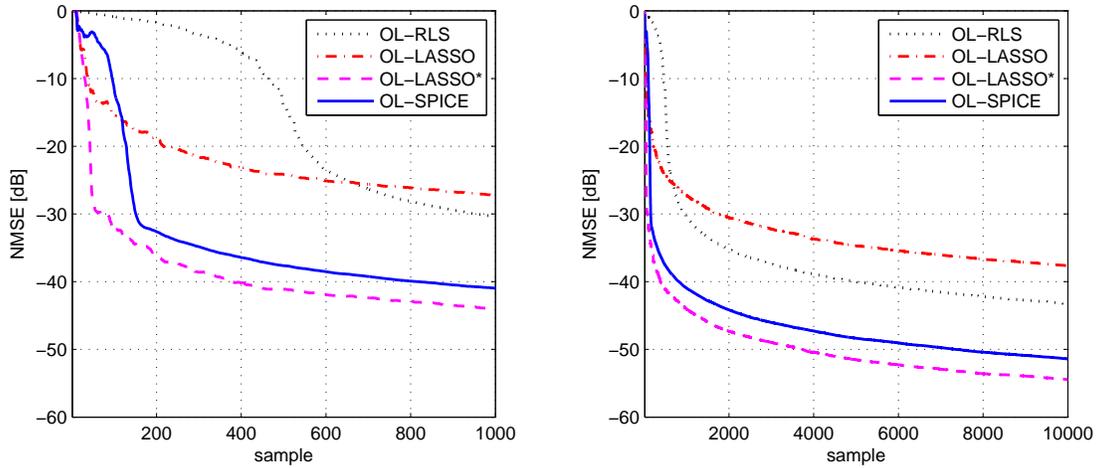}
  \end{center}
  \caption{IID regressors and deterministic $\mbs{\theta}$. NMSE versus $n$. Left: $n = 1$ to $10^3$. Right: $n = 1$ to $10^4$.  SNR=20~dB and number of nonzero parameters $p^\star = 3$. The asterisk denotes the infeasible \textsc{Ol-Lasso}.}
  \label{fig:NMSE_LASSO_iid}
\end{figure*}

We first consider $\mbs{\theta}$ to be a deterministic parameter. In the first experiment we set $p^\star = 3$ nonzero elements,
 $\theta_{10}=1$, $\theta_{20}=1$ and $\theta_{140}=3$. Note that since the regressors are drawn independently, the chosen support set for $\mbs{\theta}$ will not affect the performance. When the number of samples is very small, the estimates for \textsc{Ol-Lasso} and \textsc{Ol-Spice} have a higher variance than \textsc{Ol-Rls} which biases its estimate
more strongly towards $\mbf{0}$. For clarity we therefore set
$\hat{\mbs{\theta}}_n \equiv \mbf{0}$ during $0 \leq n \leq 20 $ for
all estimators. The results are shown in Fig.~\ref{fig:NMSE_LASSO_iid} for SNR=20~dB. We observe a significant performance gap between the feasible and
infeasible \textsc{Ol-Lasso}, which illustrates how critical it is to
tune the hyperparameter $\lambda_n$ to the unknown noise variance $\sigma^2$. Both \textsc{Ol-Lasso} and \textsc{Ol-Spice} quickly prune out many
nonzero coefficient estimates, the effect of which is visible in the transition phase of the
NMSE plot. The performance of \textsc{Ol-Rls} becomes better than that
of \textsc{Ol-Lasso} when $n > 750$. \textsc{Ol-Spice}
outperforms the feasible \textsc{Ol-Lasso} after about $n=100$ samples and is
closer to the infeasible version which uses an optimally tuned $\lambda_n$. 

\begin{figure*}
  \begin{center}
    \includegraphics[width=1.95\columnwidth]{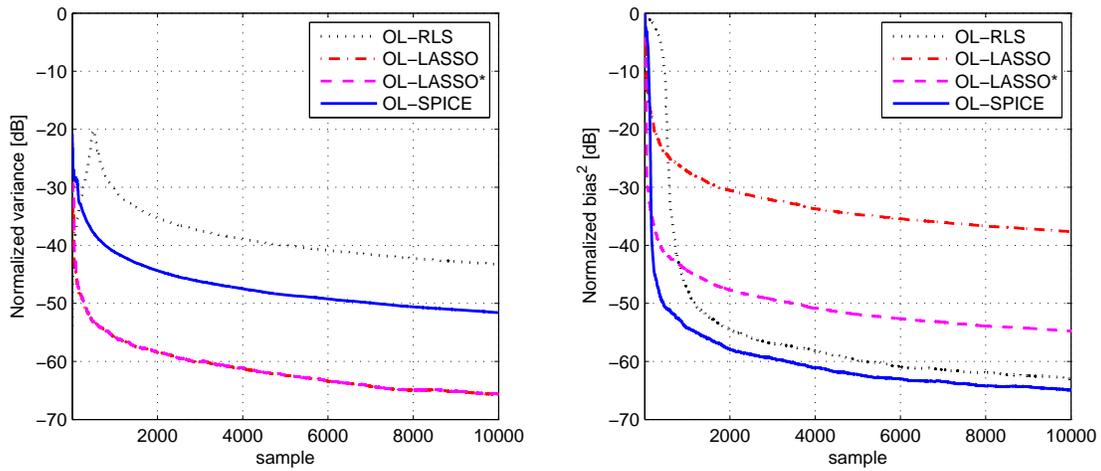}
  \end{center}
  \caption{IID regressors and deterministic $\mbs{\theta}$. Left: variance versus $n$. Right: square-bias versus $n$. SNR=20~dB and $p^\star = 3$.}
  \label{fig:NMSE_LASSO_iid_decomp}
\end{figure*}

Fig.~\ref{fig:NMSE_LASSO_iid_decomp} presents the variance and bias of the estimators by decomposing the mean square error in \eqref{eq:NMSE}. Both versions of \textsc{Ol-Lasso} exhibit much lower variance than square-bias whereas \textsc{Ol-Spice} has a more balanced variance-bias composition and noticeably the lowest bias among the considered estimators.

\begin{figure}
  \begin{center}
    \includegraphics[width=1.0\columnwidth]{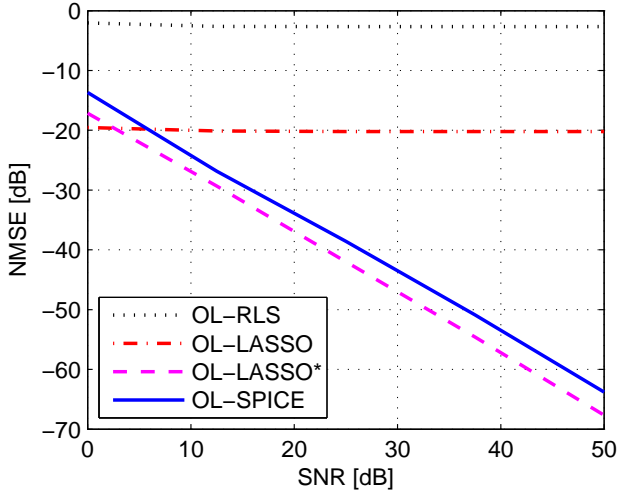}
  \end{center}
  \caption{IID regressors and deterministic $\mbs{\theta}$. NMSE versus SNR. $n=250$ samples and $p^\star = 3$. The asterisk denotes the infeasible \textsc{Ol-Lasso}.}
  \label{fig:NMSE_LASSO_SNR}
\end{figure}

In  Fig.~\ref{fig:NMSE_LASSO_SNR} we see that the NMSE for feasible \textsc{Ol-Rls} and \textsc{Ol-Lasso}, which is dominated by the bias, remains virtually unaffected by increasing SNR for a fixed number of samples $n=250$. By contrast, the errors for \textsc{Ol-Spice} and the infeasible \textsc{Ol-Lasso} decrease as the signal conditions improve.

\begin{figure}
  \begin{center}
    \includegraphics[width=1.0\columnwidth]{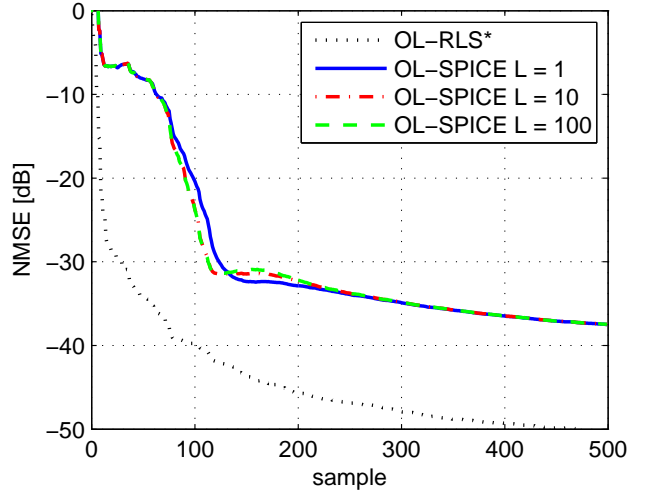}
  \end{center}
  \caption{IID regressors and deterministic $\mbs{\theta}$. NMSE versus $n$. SNR=20~dB and $p^\star = 3$. The asterisk denotes the infeasible \textsc{Ol-Rls}.}
  \label{fig:NMSE_L_iid}
\end{figure}

Next, we study the effect of the number of iteration cycles $L$ per sample
in \textsc{Ol-Spice}. The results are illustrated in
Fig.~\ref{fig:NMSE_L_iid}. We note that the performance characteristics for $L = 1$, $10$ and $100$, are very similar. For $n \leq 120$, a larger $L$ leads to slightly faster decrease of the NMSE but the differences rapidly diminish as $L$ increases and the NMSE curves almost coincide for $n
> 200$. For reference we included the infeasible \textsc{Ol-Rls}, which provides a
lower bound on the NMSE in Fig.~\ref{fig:NMSE_L_iid} and can be compared to Fig.~\ref{fig:NMSE_LASSO_iid}.

\begin{figure*}
  \begin{center}
    \includegraphics[width=1.95\columnwidth]{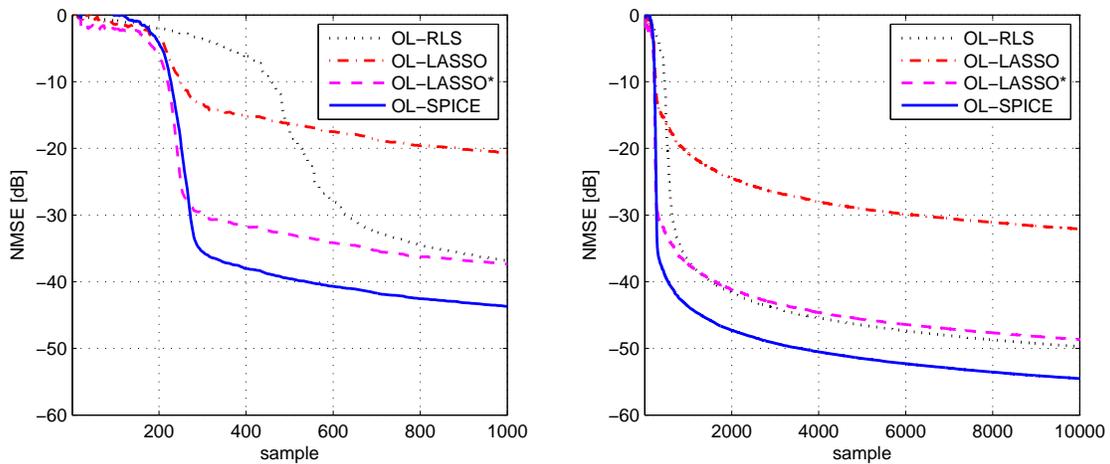}
  \end{center}
  \caption{IID regressors and deterministic $\mbs{\theta}$. NMSE versus $n$. Left: $n = 1$ to $10^3$. Right: $n = 1$ to $10^4$. SNR=20~dB and $p^\star = 50$.}
  \label{fig:NMSE_LASSO_iid_p50}
\end{figure*}

In the next experimental setup we set the number of nonzero elements
to $p^\star = 50$, $\theta_1 = \cdots = \theta_{50} = 1$, thus
increasing element density of $\mbs{\theta}$ to 10\%. The results in
Fig.~\ref{fig:NMSE_LASSO_iid_p50} show that \textsc{Ol-Spice} can
better cope with less sparse parameter vectors than
\textsc{Ol-Lasso}. For $n > 300$ it exhibits lower NMSE than the rest, owing to a lower bias (not shown here but observed by us in the numerical evaluation).

\begin{figure*}
  \begin{center}
    \includegraphics[width=1.95\columnwidth]{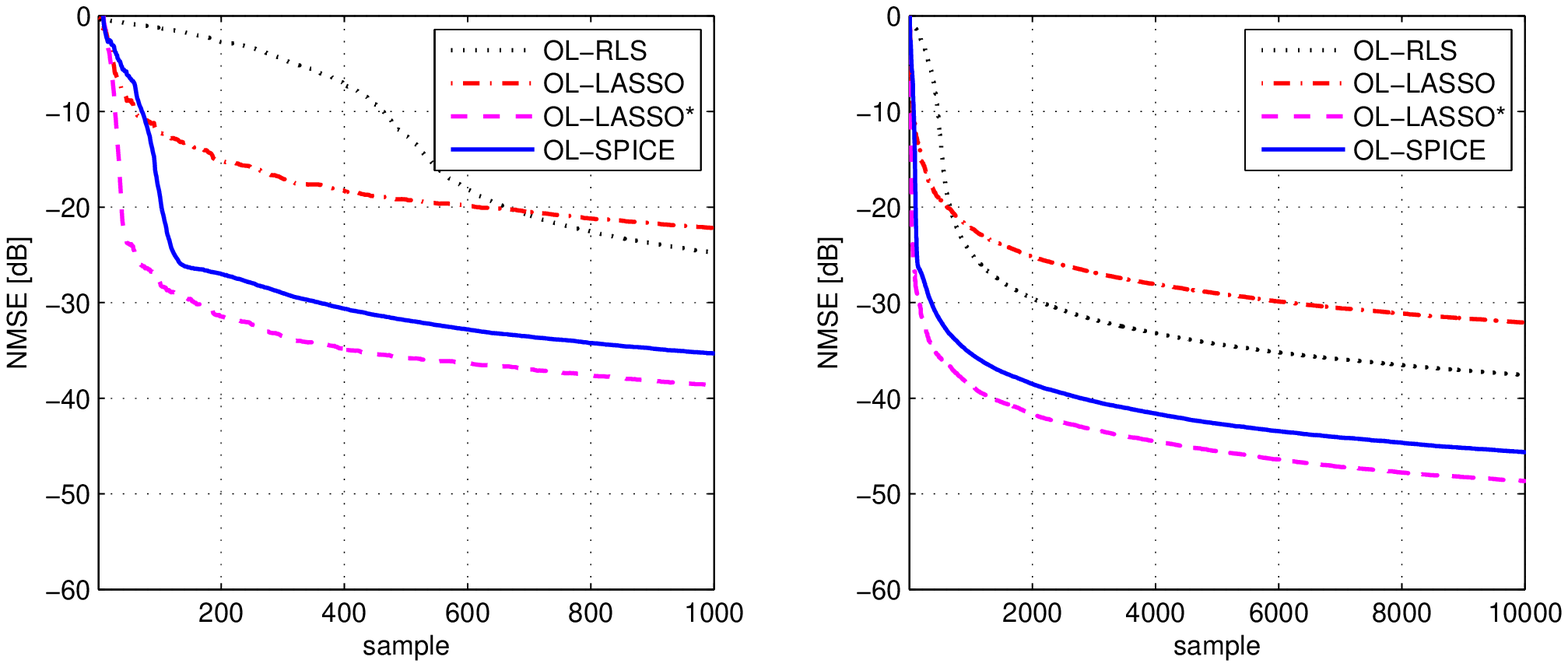}
  \end{center}
  \caption{IID regressors and stochastic $\mbs{\theta}$.  NMSE versus $n$. Left: $n = 1$ to $10^3$. Right: $n = 1$ to $10^4$. SNR=20~dB and $p^\star = 3$.}
  \label{fig:NMSE_LASSO_iid_randn}
\end{figure*}

Finally, we consider a setup where $\mbs{\theta}$ is a random parameter. Since the support set is unimportant in the present case we generate the elements $\theta_{10}$, $\theta_{20}$ and $\theta_{140}$ using independent Gaussian variables with zero-mean and unit variance, resulting in a wider dynamic range than in the previous experiments. Nevertheless, the results presented in Fig.~\ref{fig:NMSE_LASSO_iid_randn} show  performance characteristics similar to the deterministic case presented in Fig.~\ref{fig:NMSE_LASSO_iid}. 

\subsection{Real-valued example: sinusoids in noise}
In contrast to the previous example, we now present a case where the regressor columns in \eqref{eq:linearmodel}
are highly correlated. Specifically, as a further example with real-valued
parameters, we consider the identification of a
sum of $q$ sinusoids at given frequencies $\{ \omega_i \} \subset [0,\pi)$
with unknown phases $\{ \phi_i \}$ and amplitudes $\{ a_i \}$ (most of which are zero). 
In the following we will consider $q = 250$ possible sinusoids on a uniform grid of
frequencies. We set two nonzero amplitudes as $a_{10}=1$ and
$a_{20}=1$ for two slowly-varying sinusoids, narrowly spaced with $\Delta
\omega = 0.04 \pi$, and $a_{140}=3$ for a high-frequency sinusoid. The phases of the three sinusoids were set to 0.

We define the signal to noise ratio as
\begin{equation*}
\text{SNR}  = \frac{\min_{i \in \mathcal{S}} \: a^2_i}{\sigma^2},
\end{equation*}
where $\mathcal{S}$ is the set of nonzero amplitudes, 
and parameterize the signal as
\begin{equation*}
\begin{split}
y_t &= \sum^{q}_{i=1} a_i \sin(\omega_i t + \phi_i) + w_t \\
&= \sum^{q}_{i=1} A_i \cos(\omega_i t) + B_i \sin(\omega_i t)  + w_t
\\
&= \mbf{h}^\top_t \mbs{\theta} + w_t,
\end{split}
\end{equation*}
where the unknown parameter vector is $\mbs{\theta} = [A_1 \: B_1 \:
\cdots \: A_{q} \: B_{q}]^\top \in \mathbb{R}^{p}$ and $p = 2q = 500$. The
regressor vector is $\mbf{h}^\top_t = [\cos(\omega_1 t) \:
\sin(\omega_1 t) \: \cdots \: \cos(\omega_q t) \: \sin(\omega_q t)]$.

\begin{figure}
  \begin{center}
    \includegraphics[width=1.0\columnwidth]{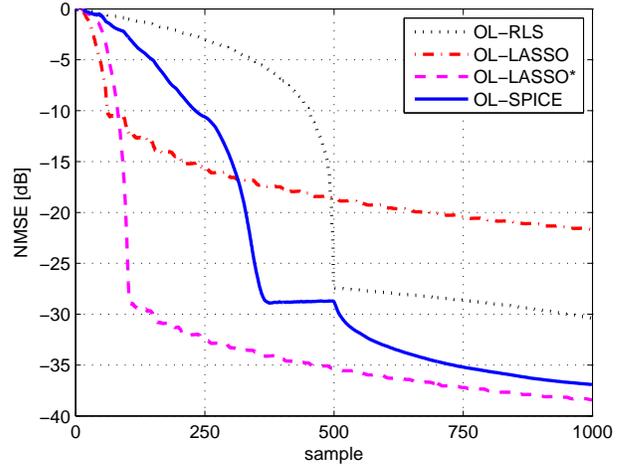}
  \end{center}
  \caption{Sinusoidal parameters $\mbs{\theta}$. NMSE versus time. SNR=20~dB. The asterisk denotes the infeasible \textsc{Ol-Lasso}.}
  \label{fig:NMSE_LASSO}
\end{figure}

\begin{figure}
  \begin{center}
    \includegraphics[width=1.0\columnwidth]{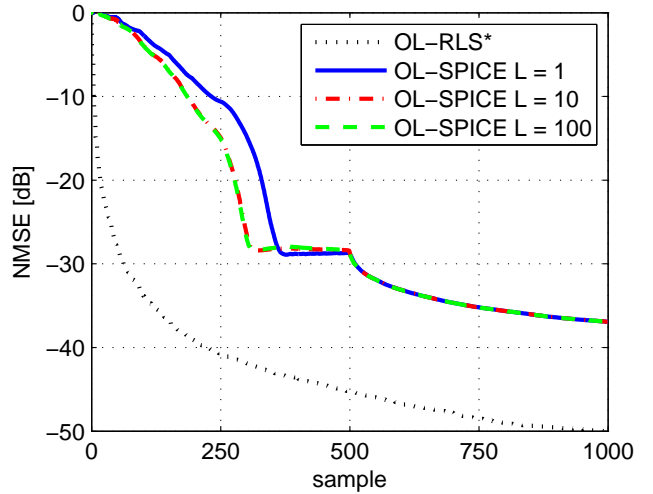}
  \end{center}
  \caption{Sinusoidal parameters $\mbs{\theta}$. NMSE versus time. SNR=20~dB. The asterisk denotes the infeasible \textsc{Ol-Rls}.}
  \label{fig:NMSE_L}
\end{figure}

We set SNR=20~dB. First, \textsc{Ol-Spice} is compared with
the feasible and infeasible \textsc{Ol-Lasso} which perform substantially
different from one another but achieve the same rate of NMSE decrease. The results are
presented in Fig.~\ref{fig:NMSE_LASSO}. For $n \leq p$,
\textsc{Ol-Spice} overtakes the feasible \textsc{Ol-Lasso} at about $n
= 300$. Notably, the NMSE of \textsc{Ol-Spice} decreases until it reaches a plateau where the
estimation errors are very small but where the noise level cannot be
properly identified. This interesting transition characteristic still awaits
a satisfactory explanation. For $n > p $, \textsc{Ol-Spice}
approaches the infeasible \textsc{Ol-Lasso} as time progresses. For reference, we
have also added the feasible \textsc{Ol-Rls} which illustrates the degradation when
not taking the parameter sparsity into account. Note, however, that \textsc{Ol-Rls} eventually
outperforms the \textsc{Ol-Lasso} estimator for which the hyperparameter has
not been finely tuned to the noise level.
Next, Fig.~\ref{fig:NMSE_L} illustrates how $L$ affects \textsc{Ol-Spice}. We see that the  performance characteristics for $L = 1$, $10$ and $100$, are very similar as was the case with weakly correlated regressor columns in Fig.~\ref{fig:NMSE_L_iid}. Setting $L=1$, however, requires slightly more samples to reach the plateau resulting in a gap in NMSE compared to $L=10$ until about $n > 375$.

\subsection{Complex-valued example: synthetic aperture radar imaging}

Finally, we illustrate how \textsc{Ol-Spice} performs in a complex-valued
case, and compare it with \textsc{Ol-Rls} and a novel form of
\textsc{Ol-Lasso} for this scenario, cf. Appendix~B. 

\begin{figure}
  \begin{center}
    \includegraphics[width=1.0\columnwidth]{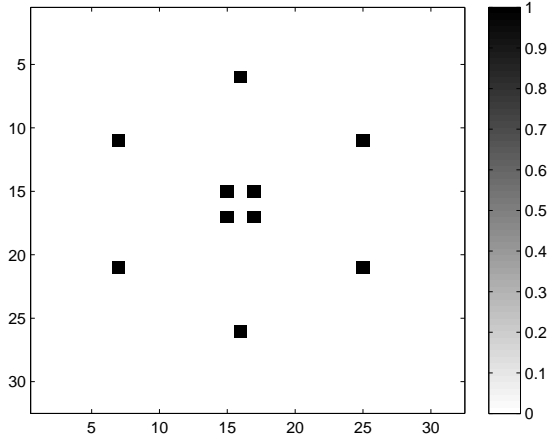}
  \end{center}
  \caption{True intensity or reflection coefficient, $\theta(\mbf{p})$, as a function of $\mbf{p} = [p_x, p_y]^\top$. 10 ideal point-scatterers are present.}
  \label{fig:SAR_reference}
\end{figure}

We consider a setup similar to that of synthetic aperture radar imaging
where an antenna transmits an electromagnetic pulse and the reflected
signal carries information about potential scatterers in the scene of
interest, cf. \cite{VuEtAl2013_bayesian}. Let $\mbf{p}$ be a position coordinate in the scene and 
$\theta(\mbf{p})$ the reflection coefficient at $\mbf{p}$. The observed 
signal is in the spatial frequency domain, where each sample
corresponds to a particular angle $\mbs{\phi}$. If we grid the
space of the scene, the signal at sample $t$ can be modeled as
\begin{equation*}
\begin{split}
y(\mbs{\phi}_t) &= \sum_{\mbf{p}} e^{-j 2 \pi
 \mbf{p}^\top \mbf{g}(\mbs{\phi}_t) } \theta(\mbf{p}) + w_t \\
&= \mbf{h}^*_t \mbs{\theta} + w_t,
\end{split} 
\end{equation*}
where $\mbs{\theta} \in \mathbb{C}^p$ is the vectorized image of
reflection coefficients. For simplicity, we
consider $\mbf{g}(\mbs{\phi}) = \mbs{\phi} \in [0,1)^2$ and
$\mbf{p} \in \mathbb{R}^2$ such that the
observation $y(\mbs{\phi}_t)$ corresponds to a coefficient of the two-dimensional
discrete Fourier transform. Here we consider the discretized scene image
to be $32 \times 32$ such that $p=1024$. The true image used in this example is shown in
Fig.~\ref{fig:SAR_reference}, which comprises $10$ point scatterers with amplitudes equal to $1$.

\begin{figure*}
  \begin{center}
   \includegraphics[width=2.0\columnwidth]{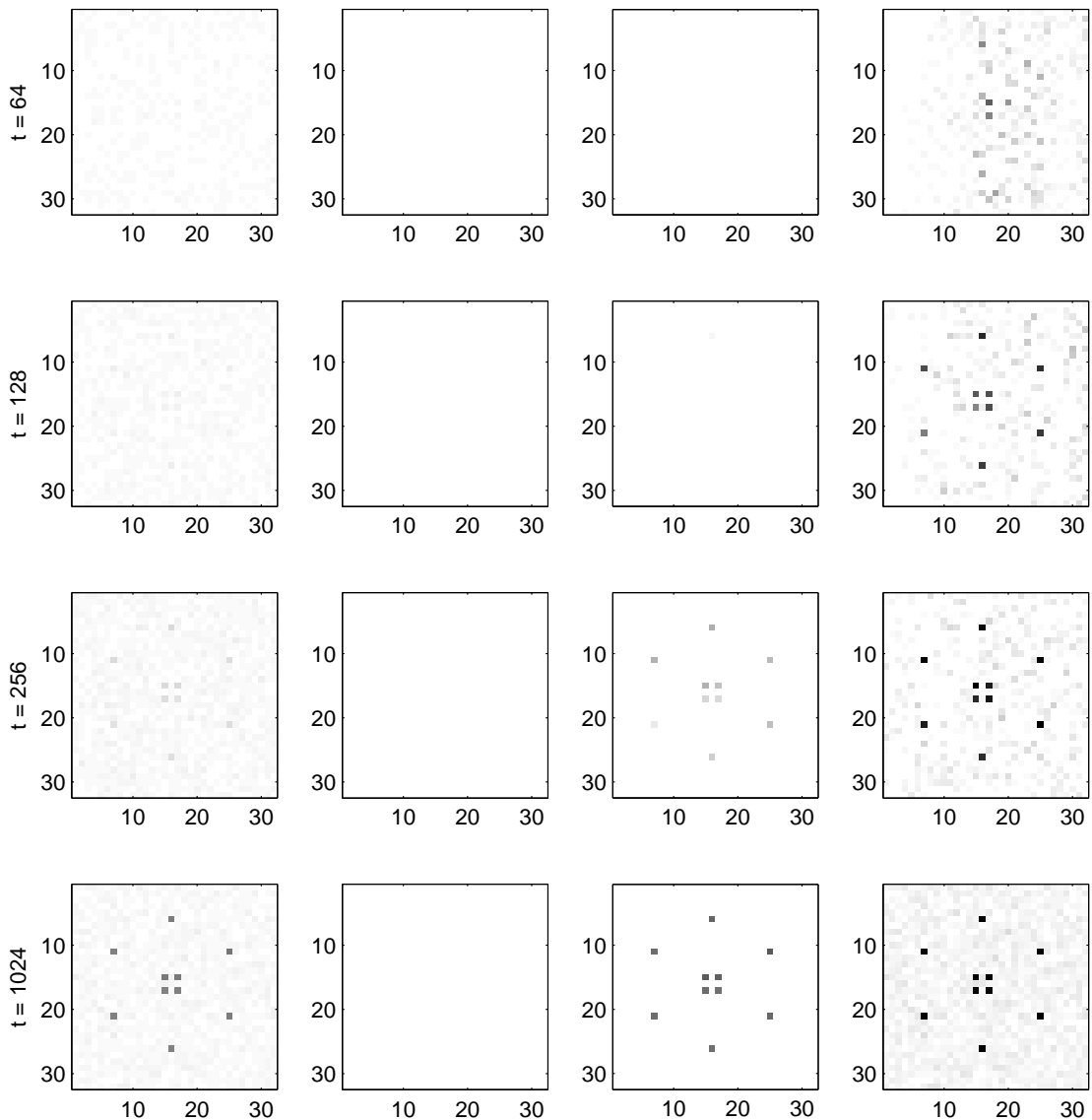}
  \end{center}
  \caption{Estimated images $|\hat{\theta}(\mbf{p})|$ at various time instants $t$ for a
    randomly chosen noise realization. The estimates for
    \textsc{Ol-Rls} and \textsc{Ol-Spice} are shown in the first and fourth columns,
    respectively. For infeasible \textsc{Ol-Lasso} with $\lambda_n =
    \sqrt{2 \sigma^2 n \log p}$ and a user-adjusted version $\lambda_n = 10^{-2}
    \sqrt{ n \log p}$, the estimates are shown in the second and third columns.  SNR=$25$~dB.}
  \label{fig:SAR_comparison}
\end{figure*}

The observations at each sample $t$ were taken at a randomly
chosen angle $\mbs{\phi}_t$ (corresponding to randomly chosen discrete spatial frequencies). The
signal to noise ratio was set to $25$~dB. In
Fig.~\ref{fig:SAR_comparison} we compare the estimated images using
\textsc{Ol-Rls}, \textsc{Ol-Lasso} and \textsc{Ol-Spice}. Note that in
this signal setup the hyperparameter in the infeasible
\textsc{Ol-Lasso} overpenalizes the $\ell_1$-norm of $\mbs{\theta}$ which
results in no visible scatterers. To produce some meaningful plots for
\textsc{Ol-Lasso}, the hyperparameter is adjusted to $\lambda_n = 10^{-2}
\sqrt{n \log p}$, which illustrates the difficulty of selecting it in
practical applications. For $t$ close to $p = 1024$, three methods estimate the locations and
intensities of the point scatterers accurately, but \textsc{Ol-Spice}
is capable of producing accurate images with far fewer samples than
the other two methods which would require fine-tuning. Indeed, the scatterer pattern is
already visible at $t=128$ samples in the \textsc{Ol-Spice} image,
without any user input.

\section{Conclusions}

We have derived an online sparse estimator, called \textsc{Ol-Spice}, that obviates the need for tuning hyperparameters. Its computational simplicity and adaptability to complex-valued parameters render it suitable for large-scale inference problems as well as real-time applications, such as system identification and synthetic aperture radar imaging. The code for \textsc{Ol-Spice} has been made available to facilitate its use in applications.


\section*{Appendix A: Linear minimum mean-square estimator and
  covariance matching}
Here we prove that the minimizer $\hat{\mbs{\theta}}$ of \eqref{eq:lmmseaugprob} is
equivalent to using the linear minimum mean-square estimator
\eqref{eq:lmmse} with covariance parameters set through covariance matching.

For notational simplicity, let $\mbs{\phi} \in \mathbb{R}^{p+1}_{++}$ denote
the covariance parameters, namely the diagonal elements of $\mbf{P}$
and $\sigma^2$, and drop subindex $n$. Further, let $\mbs{\Sigma} \triangleq \sigma^2
\mbf{I}$ so that $\mbf{R} (\mbs{\phi}) = \mbf{H} \mbf{P} \mbf{H}^* + \mbs{\Sigma}$. Now
\eqref{eq:lmmse} can be written as
\begin{equation*}
\begin{split}
\hat{\mbs{\theta}} &= \mbf{P} \mbf{H}^* \mbf{R}^{-1} \mbf{y} \\
&= (\mbf{H}^*\mbs{\Sigma}^{-1} \mbf{H} + \mbf{P}^{-1}  )^{-1} \mbf{H}^* \mbs{\Sigma}^{-1} \mbf{y}.
\end{split}
\end{equation*}
We note that \eqref{eq:lmmse} is invariant to any scaling of the
covariance parameters. That is,
\begin{equation}
\begin{split}
\hat{\mbs{\theta}}(c \mbs{\phi}) &= c \mbf{P} \mbf{H}^* (c\mbf{R})^{-1} \mbf{y} = \hat{\mbs{\theta}} (\mbs{\phi}) 
\end{split}
\label{eq:lmmseinvariance}
\end{equation}
for any $c > 0$, which follows from $\mbf{R}(c\mbs{\phi}) =  \mbf{H}_n(c\mbf{P})\mbf{H}^*+ c\mbs{\Sigma} =
c\mbf{R}(\mbs{\phi})$. Finally note that since \eqref{eq:lmmse} minimizes
\eqref{eq:lmmse_problem} it is therefore the minimizer $\hat{\mbs{\theta}}$ of the augmented problem \eqref{eq:lmmseaugprob} as well.   

We proceed by inserting \eqref{eq:lmmse} in
\eqref{eq:lmmseaugprob}; this will lead to a concentrated cost function
that is equivalent to using the covariance-matching criterion. First, using the matrix inversion lemma, note that:
\begin{equation*}
\begin{split}
\mbf{y} - \mbf{H} \hat{\mbs{\theta}} &=  \mbf{y} -
\mbf{H}(\mbf{P}^{-1} + \mbf{H}^* \mbs{\Sigma}^{-1} \mbf{H} )^{-1} \mbf{H}^* \mbs{\Sigma}^{-1}
\mbf{y} \\
&= \mbs{\Sigma} \left( \mbs{\Sigma}^{-1} - \mbs{\Sigma}^{-1} \mbf{H} (\mbf{P}^{-1 }  + \mbf{H}^* \mbs{\Sigma}^{-1} \mbf{H}
)^{-1} \mbf{H}^*\mbs{\Sigma}^{-1} \right ) \mbf{y} \\
&= \mbs{\Sigma} (\mbf{H} \mbf{P} \mbf{H}^* + \mbs{\Sigma} )^{-1} \mbf{y} \\
&= \mbs{\Sigma} \mbf{R}^{-1} \mbf{y},
\end{split}
\end{equation*}
so that
\begin{equation*}
\begin{split}
&\| \mbf{y} - \mbf{H} \hat{\mbs{\theta}} \|^2_{\mbs{\Sigma}^{-1}} + \| \hat{\mbs{\theta}}
\|^2_{\mbf{P}^{-1}} \\
&= \mbf{y}^* \mbf{R}^{-1} \mbs{\Sigma} \mbs{\Sigma}^{-1} \mbs{\Sigma}
\mbf{R}^{-1} \mbf{y} + \mbf{y}^* \mbf{R}^{-1} \mbf{H} \mbf{P} \mbf{P}^{-1} \mbf{P} \mbf{H}^* \mbf{R}^{-1} \mbf{y} \\
&= \mbf{y}^* \mbf{R}^{-1} ( \mbs{\Sigma} + \mbf{H}\mbf{P}\mbf{H}^* )\mbf{R}^{-1} \mbf{y} \\
&= \mbf{y}^* \mbf{R}^{-1} \mbf{y}.
\end{split}
\end{equation*}
Thus after concentrating out $\mbs{\theta}$, \eqref{eq:lmmseaugprob} can be
written as
\begin{equation}
\argmin_{\mbs{\phi}} \: \mbf{y}^*\mbf{R}^{-1} \mbf{y} +  \text{tr}\{ \mbf{R} \}.
\label{eq:lmmseconc}
\end{equation}

Now expand the covariance-matching criterion,
\begin{equation*}
\begin{split}
\| \mbf{R}^{-1/2}( \mbf{y}\mbf{y}^* - \mbf{R} ) \|^2_F &=  \text{tr}\{
(\mbf{y}\mbf{y}^* - \mbf{R}) \mbf{R}^{-1} (\mbf{y}\mbf{y}^* - \mbf{R}) \} \\
&= \text{tr}\{ \mbf{y}\mbf{y}^* \mbf{R}^{-1} \mbf{y}\mbf{y}^*  \} +
\text{tr}\{\mbf{R}\} - 2\text{tr}\{ \mbf{y}\mbf{y}^* \} \\
&= \mbf{y}^* \mbf{R}^{-1} \mbf{y} \| \mbf{y} \|^2 + \text{tr}\{ \mbf{R} \} + K,
\end{split}
\end{equation*}
where $K$ is a constant. The covariance matching problem can thus be
written equivalently as 
\begin{equation}
\argmin_{\mbs{\phi}} \: \mbf{y}^* \mbf{R}^{-1} \mbf{y} + \| \mbf{y} \|^{-2} \text{tr}\{ \mbf{R} \},
\label{eq:covmatch}
\end{equation}
which is similar to \eqref{eq:lmmseconc}. Let the cost
functions in \eqref{eq:lmmseconc} and \eqref{eq:covmatch} be denoted as
$J(\mbs{\phi})$ and $J'(\mbs{\phi})$, respectively. We now show that their
respective minimizers differ only by a scaling constant. That is, $\hat{\mbs{\phi}}
= c \hat{\mbs{\phi}}'$, where $c = \| \mbf{y} \|^{-1} > 0$. This follows from
\begin{equation*}
\begin{split}
c J( c \mbs{\phi} ) &= c ( \mbf{y}^*(c\mbf{R})^{-1}\mbf{y} ) + c\text{tr}\{(c\mbf{R})\} \\
&= \mbf{y}^*\mbf{R}^{-1} \mbf{y} + c^2 \text{tr}\{ \mbf{R} \} \\
&= J'(\mbs{\phi}),
\end{split}
\end{equation*}
so that for the minimizer $\hat{\mbs{\phi}}'$ we have $c J(c \hat{\mbs{\phi}}') = J'(\hat{\mbs{\phi}}') \leq
J'(\mbs{\phi}) = c J(c \mbs{\phi})$. It follows that $J(c\hat{\mbs{\phi}}') \leq
J(c\mbs{\phi})$ for all $\mbs{\phi} \in \mathbb{R}^{n+1}_{++}$, and therefore the minimizers
for the concentrated cost function  \eqref{eq:lmmseconc} and the covariance-matching
crierion \eqref{eq:covmatch} differ only by a factor $c >0$. From \eqref{eq:lmmseinvariance} we know
that the linear minimum mean-square estimator is invariant to uniform
scaling of the covariance parameters. This concludes the proof. (See
also \cite{StoicaEtAl2014_weightedspice} for other details of this result.)

\section*{Appendix B: Online Lasso for the complex-valued case}

An online cyclic \textsc{Lasso} algorithm that covers both the real and
complex-valued case can be derived using the same reparametrization employed in \textsc{Ol-Spice}. Analogous to \eqref{eq:globalcost_alt} and the derivation of \eqref{eq:J_r_phi_alt}, the cost function can be written as
\begin{equation*}
J(\theta_i) = \| \tilde{\mbf{y}}_i - \mbf{c}_i \theta_i \|^2_2 + \lambda_n |\theta_i|,
\end{equation*}
and in concentrated form,
\begin{equation*}
\begin{split}
J(r_i, \hat{\varphi}_i) &= \left( \alpha_i + \beta_i r^2_i  - 2 \gamma_i r_i \right) + \lambda_n
r_i \\
&= \beta_i \left( r^2_i - \left(\frac{2\gamma_i - \lambda_n}{\beta_i} \right) r_i \right)  + \alpha_i \\
&= \beta_i \left( r_i - \frac{1}{2}\left(\frac{2\gamma_i - \lambda_n}{\beta_i}  \right) \right)^2  + K,
\end{split}
\end{equation*}
where $K$ is a constant and the auxiliary variables can be computed as
(cf. \eqref{eq:updatevariables_compact}):
\begin{equation}
\begin{split}
\beta_i &= \Gamma^n_{ii}\\
\gamma_i &= |\zeta_i + \Gamma^n_{ii} \check{\theta}_i | .
\end{split}
\label{eq:updatevariables_compact_lasso}
\end{equation}
The minimizing argument $r_i \geq 0$ is given by
\begin{equation}
r_i = \text{max}\left( \frac{2\gamma_i - \lambda_n}{2 \beta_i}, 0 \right).
\label{eq:r_i_gen_lasso}
\end{equation}
Thus we have the minimizer $\hat{\theta}_i = \hat{r}_i e^{j
  \hat{\varphi}_i}$, where
\begin{equation}
\hat{\varphi}_i = \text{arg}(\zeta_i + \Gamma^n_{ii} \check{\theta}_i).
\label{eq:phi_hat_update_lasso}
\end{equation}
Note that the above derivation does not involve the variable
$\alpha_i$ from Section~\ref{sec:onlinespice} or the variables $\eta_n$ and $\kappa^n$ in
the online formulation of \textsc{Spice}, cf. \eqref{eq:updatevariables_z},
\eqref{eq:auxiliaryvariables}, \eqref{eq:recursivevariables} and
\eqref{eq:updatevariables_compact}. The result is summarized in Algorithm~\ref{alg:onlineslasso}.
\begin{algorithm}
  \caption{Online \textsc{Lasso}} \label{alg:onlineslasso}
\begin{algorithmic}[1]
    \State Input: $y_n$, $\mbf{h}_n$ and $\check{\mbs{\theta}}$
    \State $\mbs{\Gamma} := \mbs{\Gamma} + \mbf{h}_n \mbf{h}^*_n$
    \State $\mbs{\rho} := \mbs{\rho} + \mbf{h}_n y_n$ 
  \State $\mbs{\zeta} = \mbs{\rho} - \mbs{\Gamma} \check{\mbs{\theta}} $
    \Repeat
        \State $i = 1, \dots, p$
        \State Compute \eqref{eq:updatevariables_compact_lasso},
        \eqref{eq:r_i_gen_lasso} and \eqref{eq:phi_hat_update_lasso}
        \State $\hat{\theta}_i = \hat{r}_i e^{j
  \hat{\varphi}_i}$
       \State  $\mbs{\zeta}  := \mbs{\zeta} + [\mbs{\Gamma}]_i (\check{\theta}_i - \hat{\theta}_i)$ 
        \State $\check{\theta}_i := \hat{\theta}_i$
    \Until{ termination }
    \State Output: $\hat{\mbs{\theta}}$
\end{algorithmic}
\end{algorithm}

\begin{IEEEbiography}{Dave Zachariah} is a researcher at Uppsala University in Sweden. He received the M.S. degree in electrical engineering from Royal Institute of Technology (KTH), Stockholm, Sweden, in 2007. During 2007-8 he worked as a research engineer at Global IP Solutions, Stockholm. He received the Tech. Lic. and Ph.D. degrees in signal processing from KTH in 2011 and 2013, respectively.

His research interests include statistical signal processing, machine learning, sensor
fusion, and localization.
\end{IEEEbiography}

\begin{IEEEbiography}{Petre Stoica} is a researcher and educator in the field of signal processing and its applications to radar/sonar, communications and biomedicine. He is a
professor of signal and system modeling at Uppsala University in Sweden, and a member of the Royal Swedish Academy of Engineering Sciences, the Romanian Academy (honorary), the European Academy of Sciences, and the Royal Society of Sciences in Uppsala.
\end{IEEEbiography}

\bibliographystyle{ieeetr}
\bibliography{refs_onlinesparse}

\begin{thebibliography}{10}

\bibitem{BourguignonEtAl2007_sparsity}
S.~Bourguignon, H.~Carfantan, and J.~Idier, ``A sparsity-based method for the
  estimation of spectral lines from irregularly sampled data,'' {\em IEEE J.
  Selected Topics in Signal Processing}, vol.~1, no.~4, pp.~575--585, 2007.

\bibitem{KleinEtAl2008_sparsity}
A.~Klein, H.~Carfantan, D.~Testa, A.~Fasoli, and J.~Snipes, ``A sparsity-based
  method for the analysis of magnetic fluctuations in unevenly-spaced {M}irnov
  coils,'' {\em Plasma Physics and Controlled Fusion}, vol.~50, no.~12,
  p.~125005, 2008.

\bibitem{StoicaEtAl2011_spicespectral}
P.~Stoica, P.~Babu, and J.~Li, ``New method of sparse parameter estimation in
  separable models and its use for spectral analysis of irregularly sampled
  data,'' {\em IEEE Trans. Signal Processing}, vol.~59, no.~1, pp.~35--47,
  2011.

\bibitem{Stoica&Babu2012_spice}
P.~Stoica and P.~Babu, ``{SPICE} and {LIKES}: Two hyperparameter-free methods
  for sparse-parameter estimation,'' {\em Signal Processing}, vol.~92, no.~7,
  pp.~1580--1590, 2012.

\bibitem{Gorodnitsky&Rao1997_sparse}
I.~F. Gorodnitsky and B.~D. Rao, ``Sparse signal reconstruction from limited
  data using {FOCUSS}: A re-weighted minimum norm algorithm,'' {\em IEEE Signal
  Processing}, vol.~45, no.~3, pp.~600--616, 1997.

\bibitem{MalioutovEtAl2005_sparse}
D.~Malioutov, M.~{\c{C}}etin, and A.~S. Willsky, ``A sparse signal
  reconstruction perspective for source localization with sensor arrays,'' {\em
  IEEE Trans. Signal Processing}, vol.~53, no.~8, pp.~3010--3022, 2005.

\bibitem{StoicaEtAl2011_spicearray}
P.~Stoica, P.~Babu, and J.~Li, ``{SPICE}: A sparse covariance-based estimation
  method for array processing,'' {\em IEEE Trans. Signal Processing}, vol.~59,
  no.~2, pp.~629--638, 2011.

\bibitem{Tibshirani1996_lasso}
R.~Tibshirani, ``Regression shrinkage and selection via the lasso,'' {\em
  Journal of the Royal Statistical Society. Series B (Methodological)},
  pp.~267--288, 1996.

\bibitem{WuEtAl2009_genome}
T.~T. Wu, Y.~F. Chen, T.~Hastie, E.~Sobel, and K.~Lange, ``Genome-wide
  association analysis by lasso penalized logistic regression,'' {\em
  Bioinformatics}, vol.~25, no.~6, pp.~714--721, 2009.

\bibitem{LuEtAl2011_lasso}
Y.~Lu, Y.~Zhou, W.~Qu, M.~Deng, and C.~Zhang, ``A lasso regression model for
  the construction of micro{RNA}-target regulatory networks,'' {\em
  Bioinformatics}, vol.~27, no.~17, pp.~2406--2413, 2011.

\bibitem{LustigEtAl2008_compressedmri}
M.~Lustig, D.~L. Donoho, J.~M. Santos, and J.~M. Pauly, ``Compressed sensing
  {MRI},'' {\em IEEE Signal Processing Magazine}, vol.~25, no.~2, pp.~72--82,
  2008.

\bibitem{DonevaEtAl2010_compressed_mri}
M.~Doneva, P.~B{\"o}rnert, H.~Eggers, C.~Stehning, J.~S{\'e}n{\'e}gas, and
  A.~Mertins, ``Compressed sensing reconstruction for magnetic resonance
  parameter mapping,'' {\em Magnetic Resonance in Medicine}, vol.~64, no.~4,
  pp.~1114--1120, 2010.

\bibitem{ChenEtAl2009_sparse}
Y.~Chen, Y.~Gu, and A.~O. Hero, ``Sparse {LMS} for system identification,'' in
  {\em IEEE Int. Conf. Acoustics, Speech and Signal Processing (ICASSP)},
  pp.~3125--3128, IEEE, 2009.

\bibitem{KalouptsidisEtAl2011_adaptive}
N.~Kalouptsidis, G.~Mileounis, B.~Babadi, and V.~Tarokh, ``Adaptive algorithms
  for sparse system identification,'' {\em Signal Processing}, vol.~91, no.~8,
  pp.~1910--1919, 2011.

\bibitem{KopsinisEtAl2011_online}
Y.~Kopsinis, K.~Slavakis, and S.~Theodoridis, ``Online sparse system
  identification and signal reconstruction using projections onto weighted
  balls,'' {\em IEEE Trans. Signal Processing}, vol.~59, no.~3, pp.~936--952,
  2011.

\bibitem{Glentis2014_adaptiveslim}
G.-O. Glentis, ``Adaptive identification of sparse systems using the {SLIM},''
  in {\em Signal Processing Conference (EUSIPCO), 2014 Proceedings of the 22nd
  European}.

\bibitem{ThemelisEtAl2014_variationalbayessparse}
K.~Themelis, A.~Rontogiannis, and K.~Koutroumbas, ``A variational {B}ayes
  framework for sparse adaptive estimation,'' {\em Signal Processing, IEEE
  Transactions on}, vol.~62, pp.~4723--4736, Sept 2014.

\bibitem{Zhu&Bamler2010_tomographic}
X.~X. Zhu and R.~Bamler, ``Tomographic {SAR} inversion by l1-norm
  regularization---the compressive sensing approach,'' {\em IEEE Trans.
  Geoscience and Remote Sensing}, vol.~48, no.~10, pp.~3839--3846, 2010.

\bibitem{CetinEtAl2014_sparsitysar}
M.~Cetin, I.~Stojanovic, O.~Onhon, K.~Varshney, S.~Samadi, W.~Karl, and
  A.~Willsky, ``Sparsity-driven synthetic aperture radar imaging:
  Reconstruction, autofocusing, moving targets, and compressed sensing,'' {\em
  IEEE Signal Processing Magazine}, vol.~31, no.~4, pp.~27--40, 2014.

\bibitem{Malecki&Donoho2010_tunedcompressed}
A.~Maleki and D.~Donoho, ``Optimally tuned iterative reconstruction algorithms
  for compressed sensing,'' {\em IEEE J. Selected Topics in Signal Processing},
  vol.~4, no.~2, pp.~330--341, 2010.

\bibitem{GiraudEtAl_2012_unknownvariance}
C.~Giraud, S.~Huet, and N.~Verzelen, ``High-dimensional regression with unknown
  variance,'' {\em Statistical Science}, vol.~27, no.~4, pp.~500--518, 2012.

\bibitem{YangEtAl2012_compressedsensingcomplex}
Z.~Yang, C.~Zhang, and L.~Xie, ``On phase transition of compressed sensing in
  the complex domain,'' {\em IEEE Signal Processing Letters}, vol.~19,
  pp.~47--50, Jan 2012.

\bibitem{MalekiEtAl2013_complexlasso}
A.~Maleki, L.~Anitori, Z.~Yang, and R.~Baraniuk, ``Asymptotic analysis of
  complex {LASSO} via complex approximate message passing {(CAMP)},'' {\em IEEE
  Trans. Information Theory}, vol.~59, pp.~4290--4308, July 2013.

\bibitem{Bjorck1996_numericalls}
A.~Bj{\"o}rck, {\em Numerical methods for least squares problems}.
\newblock Siam, 1996.

\bibitem{Soderstrom&Stoica1988_system}
T.~S{\"o}derstr{\"o}m and P.~Stoica, {\em System identification}.
\newblock Prentice-Hall, Inc., 1988.

\bibitem{KailathEtAl2000_linear}
T.~Kailath, A.~H. Sayed, and B.~Hassibi, {\em Linear estimation}.
\newblock Prentice-Hall, Inc., 2000.

\bibitem{Stoica&Ahgren2002_exactrls}
P.~Stoica and P.~{\AA}hgren, ``Exact initialization of the recursive
  least-squares algorithm,'' {\em International Journal of Adaptive Control and
  Signal Processing}, vol.~16, no.~3, pp.~219--230, 2002.

\bibitem{DumitrescuEtAl2012_greedyrls}
B.~Dumitrescu, A.~Onose, P.~Helin, and I.~Tabus, ``Greedy sparse {RLS},'' {\em
  IEEE Trans. Signal Processing}, vol.~60, no.~5, pp.~2194--2207, 2012.

\bibitem{Eksioglu&Korhan2011_rlsconvex}
E.~M. Eksioglu and A.~K. Tanc, ``{RLS} algorithm with convex regularization,''
  {\em IEEE Signal Processing Letters}, vol.~18, no.~8, pp.~470--473, 2011.

\bibitem{ChenEtAl1998_bpdn}
S.~S. Chen, D.~L. Donoho, and M.~A. Saunders, ``Atomic decomposition by basis
  pursuit,'' {\em SIAM journal on scientific computing}, vol.~20, no.~1,
  pp.~33--61, 1998.

\bibitem{Fu1998_penalized}
W.~J. Fu, ``Penalized regressions: the bridge versus the lasso,'' {\em Journal
  of computational and graphical statistics}, vol.~7, no.~3, pp.~397--416,
  1998.

\bibitem{FriedmanEtAl2007_pathwise}
J.~Friedman, T.~Hastie, H.~H{\"o}fling, and R.~Tibshirani, ``Pathwise
  coordinate optimization,'' {\em The Annals of Applied Statistics}, vol.~1,
  no.~2, pp.~302--332, 2007.

\bibitem{BabadiEtAl2010_sparls}
B.~Babadi, N.~Kalouptsidis, and V.~Tarokh, ``{SPARLS}: The sparse {RLS}
  algorithm,'' {\em IEEE Trans. Signal Processing}, vol.~58, no.~8,
  pp.~4013--4025, 2010.

\bibitem{Garrigues&Ghaoi2009_homotopylasso}
P.~Garrigues and L.~E. Ghaoui, ``An homotopy algorithm for the lasso with
  online observations,'' in {\em Advances in neural information processing
  systems}, pp.~489--496, 2009.

\bibitem{Asif&Romberg2010_homotopylasso}
M.~Salman~Asif and J.~Romberg, ``Dynamic updating for $\ell_{1}$
  minimization,'' {\em IEEE J. Selected Topics in Signal Processing}, vol.~4,
  pp.~421--434, April 2010.

\bibitem{AngelosanteEtAl2010_onlinelasso}
D.~Angelosante, J.~A. Bazerque, and G.~B. Giannakis, ``Online adaptive
  estimation of sparse signals: Where {RLS} meets the $\ell_1$-norm,'' {\em
  IEEE Trans. Signal Processing}, vol.~58, no.~7, pp.~3436--3447, 2010.

\bibitem{Fuchs2005_recovery}
J.-J. Fuchs, ``Recovery of exact sparse representations in the presence of
  bounded noise,'' {\em IEEE Trans. Information Theory}, vol.~51, no.~10,
  pp.~3601--3608, 2005.

\bibitem{DonohoEtAl2006_stable}
D.~L. Donoho, M.~Elad, and V.~N. Temlyakov, ``Stable recovery of sparse
  overcomplete representations in the presence of noise,'' {\em IEEE Trans.
  Information Theory}, vol.~52, no.~1, pp.~6--18, 2006.

\bibitem{Tropp2006_justrelax}
J.~A. Tropp, ``Just relax: Convex programming methods for identifying sparse
  signals in noise,'' {\em IEEE Trans. Information Theory}, vol.~52, no.~3,
  pp.~1030--1051, 2006.

\bibitem{BelloniEtAl2011_squarerootlasso}
A.~Belloni, V.~Chernozhukov, and L.~Wang, ``Square-root lasso: pivotal recovery
  of sparse signals via conic programming,'' {\em Biometrika}, vol.~98, no.~4,
  pp.~791--806, 2011.

\bibitem{Sun&Zhang2012_scaledsparse}
T.~Sun and C.-H. Zhang, ``Scaled sparse linear regression,'' {\em Biometrika},
  p.~ass043, 2012.

\bibitem{Huber2011_robust}
P.~J. Huber, {\em Robust statistics}.
\newblock Springer, 2011 [1981].

\bibitem{vanTrees2013_detection}
H.~Van~Trees and K.~Bell, {\em Detection Estimation and Modulation Theory,
  Pt.I}.
\newblock Detection Estimation and Modulation Theory, Wiley, second~ed., 2013
  [1968].

\bibitem{Anderson1989_linear}
T.~W. Anderson, ``Linear latent variable models and covariance structures,''
  {\em Journal of Econometrics}, vol.~41, no.~1, pp.~91--119, 1989.

\bibitem{OtterstenEtAl1998_covariance}
B.~Ottersten, P.~Stoica, and R.~Roy, ``Covariance matching estimation
  techniques for array signal processing applications,'' {\em Digital Signal
  Processing}, vol.~8, no.~3, pp.~185--210, 1998.

\bibitem{RojasEtAl2013_spicenote}
C.~Rojas, D.~Katselis, and H.~Hjalmarsson, ``A note on the {SPICE} method,''
  {\em IEEE Trans. Signal Processing}, vol.~61, pp.~4545--4551, Sept 2013.

\bibitem{Babu&Stoica2014_connection}
P.~Babu and P.~Stoica, ``Connection between {SPICE} and square-root {LASSO} for
  sparse parameter estimation,'' {\em Signal Processing}, vol.~95, pp.~10--14,
  2014.

\bibitem{StoicaEtAl2014_weightedspice}
P.~Stoica, D.~Zachariah, and J.~Li, ``Weighted {SPICE}: A unifying approach for
  hyperparameter-free sparse estimation,'' {\em Digital Signal Processing},
  vol.~33, pp.~1--12, 2014.

\bibitem{VuEtAl2013_bayesian}
D.~Vu, M.~Xue, X.~Tan, and J.~Li, ``A {B}ayesian approach to {SAR} imaging,''
  {\em Digital Signal Processing}, vol.~23, no.~3, pp.~852--858, 2013.

\end{thebibliography}

\end{document}